\documentclass[11pt,a4paper]{amsart}
\usepackage{amssymb,amsmath,amsthm}
\newtheorem{theorem}{Theorem}[section]
\newtheorem{lemma}[theorem]{Lemma}
\newtheorem{cor}[theorem]{Corollary}
\newtheorem{prop}[theorem]{Proposition}
\newtheorem{scholium}[theorem]{Scholium}
\theoremstyle{remark}
\newtheorem{remark}[theorem]{Remark}
\newtheorem{example}[theorem]{Example}
\theoremstyle{definition}
\newtheorem{definition}[theorem]{Definition}
\newtheorem{notation}[theorem]{Notation}
\numberwithin{equation}{section}
\DeclareMathOperator{\Aut}{Aut}
\DeclareMathOperator{\Ad}{Ad}
\DeclareMathOperator{\clsp}{\overline{span}}
\DeclareMathOperator{\ran}{ran}
\DeclareMathOperator{\tail}{Tail}
\DeclareMathOperator{\fpa}{FPA}
\newcommand{\rep}[1]{\tilde\pi[#1]}
\newcommand{\abs}[1]{\lvert#1\rvert}
\newcommand{\norm}[1]{\lVert#1\rVert}
\newcommand{\cstar}{$C^*$\ndash}
\newcommand{\Star}{${}^*$\ndash}
\newcommand{\ue}{\stackrel{u}{\sim}}
\newcommand{\qe}{\stackrel{q}{\sim}}
\newcommand{\conj}{\stackrel{c}{\sim}}
\newcommand{\set}[1]{\{#1\}}
\newcommand{\setspace}[1]{\{\,#1\,\}}
\newcommand{\ip}[2]{\langle #1, #2 \rangle}
\newcommand{\ndash}{\nobreakdash-}
\newcommand{\Ra}{\Rightarrow}
\newcommand{\bh}{\Bb(\Hh)}
\newcommand{\F}{{\Ff}_n}
\newcommand{\On}{{\Oo}_n}
\newcommand{\field}[1]{\mathbb{#1}}
\newcommand{\CC}{\field{C}}
\newcommand{\NN}{\field{N}}
\newcommand{\TT}{\field{T}}
\newcommand{\ZZ}{\field{Z}}
\newcommand{\Bb}{{\mathcal B}}
\newcommand{\Dd}{{\mathcal D}}
\newcommand{\Ee}{{\mathcal E}}
\newcommand{\Ff}{{\mathcal F}}
\newcommand{\Gg}{{\mathcal G}}
\newcommand{\Hh}{{\mathcal H}}
\newcommand{\Ii}{{\mathcal I}}
\newcommand{\Jj}{{\mathcal J}}
\newcommand{\Kk}{{\mathcal K}}
\newcommand{\Mm}{{\mathcal M}}
\newcommand{\Oo}{{\mathcal O}}
\newcommand{\Ss}{{\mathcal S}}
\newcommand{\Uu}{{\mathcal U}}

\begin{document}
%
%
\title[endomorphisms of $\bh$]{ Endomorphisms of $\bh$,
extensions of pure states, and a class of representations of 
$\mathcal O_n$}
\author{Neal J. Fowler}
\address{Department of Mathematics  \\
      University of Newcastle\\  Callaghan, NSW  2308 \\ AUSTRALIA}
\email{neal@maths.newcastle.edu.au}
\author{Marcelo Laca}
\email{marcelo@math.newcastle.edu.au}
\date{24 September 1997}
\begin{abstract}
Let $\F$ be the fixed-point algebra of the gauge 
action of the circle on the Cuntz algebra $\On$.
For every pure state $\rho$ of $\F$ and every representation $\theta$
of $C(\TT)$ we construct a representation of $\On$, and we use
the resulting class of representations to parameterize the space of all
states of $\On$ which extend $\rho$.
We show that the gauge group acts transitively on
the pure extensions of $\rho$ and that the action is $p$-to-$1$
with $p$ the period of $\rho$ under the usual shift.
We then use the above representations of $\On$ to construct endomorphisms
of $\bh$, which we classify up to conjugacy
in terms of the parameters $\rho$ and $\theta$.
In particular our construction yields every ergodic endomorphism $\alpha$
whose tail algebra $\bigcap_k \alpha^k(\bh)$ has a minimal projection,
and our results classify these ergodic endomorphisms by
an equivalence relation on the pure states of $\F$. 
As examples we analyze the ergodic endomorphisms arising from periodic
pure product states of $\F$, for which we are able to give 
a geometric complete conjugacy invariant, generalizing results of
Stacey \cite{sta-shift}, Laca \cite{laca1,laca2}, 
and Bratteli-J{\o}rgensen-Price \cite{bjp} on the shifts of Powers \cite{powers}.
\end{abstract}
\maketitle
%
%
\section*{Introduction.}
Let $\bh$ denote the algebra of bounded linear operators
on a separable complex Hilbert space $\Hh$.
An {\em endomorphism\/} of $\bh$ is a homomorphism of
$\bh$ into itself which preserves adjoints.
The main goal of this paper is to find complete conjugacy
invariants for a certain class of endomorphisms of $\bh$.
Two endomorphisms are {\em conjugate\/} if there is an intertwining
isomorphism of the underlying algebras.  Every
isomorphism of $\bh$ is unitarily implemented, so conjugacy
for endomorphisms of $\bh$ is {\em spatial\/} equivalence,
the strongest reasonable equivalence relation in any classification scheme.

If an endomorphism $\alpha$ fixes only the scalar operators
it is called {\em ergodic\/}, and if its {\em tail algebra\/}
$\bigcap_k\alpha^k(\bh)$ consists only of scalars it is called
a {\em shift\/}.
The class of endomorphisms we shall consider includes 
every ergodic endomorphism whose tail algebra has minimal projection;
in particular, it includes all shifts.

We shall only consider endomorphisms which preserve
the identity operator $I$ on $\Hh$.  At the other extreme
are those endomorphisms $\alpha$ which are completely
nonunital in the sense that $\alpha^k(I)$ decreases strongly
to zero; we refer the reader to \cite[\S2]{laca1} for the
classification of such endomorphisms.  Since any endomorphism
can be decomposed into unital and completely nonunital components
which determine its conjugacy class, our focus on unital endomorphisms
is justified.

As is customary, for $2 \le n \le \infty$ we denote by $\On$
the \cstar algebra defined by Cuntz in \cite{cun}.  Let
$\setspace{v_a: 1\le a \le n}$ denote the distinguished generating
isometries in $\On$, so that $\sum_{a=1}^n v_av_a^* \le 1$,
with equality if $n$ is finite;
$\On$ is the universal \cstar algebra generated by such collections
of isometries.  There is a correspondence between
endomorphisms of $\bh$ and representations of Cuntz algebras, 
which stems from the observation by Arveson \cite{arv}
that every endomorphism $\alpha$ of $\bh$ can be implemented by
a collection $S_1$, \dots, $S_n$ of isometries on $\Hh$
via $\alpha(A) = \sum S_aAS_a^*$;
if $\alpha$ is unital then such a collection gives rise
to a representation of $\On$ via $v_a \mapsto S_a$.
Conversely, any representation $\pi:\On\to\bh$ gives rise
to an endomorphism $\Ad\pi$ of $\bh$ via
\begin{equation}\label{eq:adpi}
\Ad\pi(A) := \sum_{a=1}^n \pi(v_a)A\pi(v_a)^*, \qquad A\in\bh; \tag{$\dag$}
\end{equation}
for infinite $n$ the above sum converges in the strong operator topology
for every $A$.
If $n<\infty$ then $\Ad\pi$ is unital, but for infinite
$n$ this need not be the case.  A representation $\pi$
of $\Oo_\infty$ for which $\Ad\pi$ is unital (i.e., for which
$\sum_{a=1}^\infty \pi(v_a)\pi(v_a)^*$ converges strongly to $I$)
is called {\em essential\/}.  

There is an obvious way to construct endomorphisms from states
of $\On$: use the GNS representation for the state to implement
an endomorphism via \eqref{eq:adpi}. 
This correspondence allows us to study endomorphisms by looking
at states of $\On$; e.g. ergodic endomorphisms arise from pure states,
and conjugacy of endomorphisms corresponds to quasi-free equivalence
of states \cite{laca1}.

A commonly used method
of analyzing $\On$ is to exploit the {\em gauge action\/}
$\gamma$ of the circle $\TT$ on $\On$ determined by
$\gamma_\lambda(v_a) = \lambda v_a$.  We will denote
by $\F$ the fixed-point algebra of this action.
When $n$ is finite, $\F$ is canonically isomorphic to
the UHF algebra $M_n\otimes M_n \otimes M_n \otimes\dotsb$, and hence
carries a canonical unital shift, given at the \cstar algebra level
by a formula analogous to \eqref{eq:adpi}.  This shift does not exist on
$\Ff_\infty$ because the strong sum does not make sense at this level, 
but one can always shift a state $\rho$ of $\F$
(for finite or infinite $n$) by defining 
\[
\alpha^*\rho(x) := \sum_{a=1}^n \rho(v_axv_a^*),\qquad x\in\F.
\]
A state $\tilde\rho$ of $\Oo_\infty$ which extends $\rho$ is essential
(i.e., its GNS representation is essential) if and only if each of the shifted
positive linear functionals $\alpha^{*k}\rho$ is a state \cite[Remark~2.10]{laca1};
in this case we say that $\rho$ is {\em essential\/}.  We will only consider
essential states of $\F$, with the understanding that the adjective
is superfluous when $n$ is finite.

The state space of $\F$ is more tractable than that of $\On$,
and has often been used to study representations of Cuntz algebras
\cite{evans,ace,spi-gi} and unital endomorphisms of $\bh$ \cite{laca1,laca2,bjp,bj}.
Having a specific procedure for extending states of $\F$ to $\On$
is extremely useful, especially if it allows one to apply Powers'
criteria for states of UHF algebras \cite{pow-th} to decide when
an extension is pure and when different extensions are unitarily equivalent.

Perhaps the most obvious way to extend a state is by
composition with the canonical conditional expectation $\Phi : \On \to \F$
obtained by averaging over the gauge action.
This gives the unique gauge-invariant state
of $\On$ which extends the given state of $\F$. 
Gauge-invariant extensions of states of $\F$ have been 
considered before: e.g., by Evans \cite{evans}, by
Araki, Carey and Evans for product states and
$n < \infty$ \cite{ace}, and later by Laca for factor states
and any $n$ \cite{laca2}.  Extensions of diagonal states
(i.e., states of the diagonal subalgebra $\Dd$ of $\F$) 
have been considered by Spielberg \cite{spi-gi,spi-qu},
by Archbold, Lazar, Tsui and Wright \cite{altw},
and by  Stacey \cite{sta-gi} 
in the context of extending the trace on the Choi subalgebra of $\Oo_2$.
In earlier work, Lazar, Tsui and Wright \cite{ltw} dealt with pure state extensions 
of pure diagonal states, and identified the unique pure state extension
of a nonperiodic (irrational) point in the spectrum of $\Dd$ \cite{cun}. 

A different procedure for extending pure states of $\F$ was given in
Theorem~4.3 of \cite{laca1}, where it played a key r\^ole in
the classification of shifts of $\bh$ up to conjugacy. 
Roughly, the technique used there consisted of lifting the GNS
representation of a state from $\F$ to $\On$ 
without changing the Hilbert space; see also \cite[\S5]{bjp}.
These two techniques for extending a pure state $\rho$ of $\F$
are in a sense opposite:  the first one always works, but only
gives an extension which is pure if $\rho$ is {\em aperiodic\/}
in the sense that its translates by powers of the canonical shift $\alpha^*$
are mutually disjoint; the second only works if $\rho$ is
{\em quasi-invariant\/} in the sense that it is quasi-equivalent
to $\alpha^*\rho$, but then gives extensions which are pure.

Here we show how to extend any pure state $\rho$ of $\F$ which is
{\em periodic\/} in the sense that $\rho$ is quasi-equivalent to $\alpha^{*p}\rho$
for some positive integer $p$. Our procedure
interpolates between the two techniques described above,
and explains them as extreme cases of the same construction.

The paper is organized as follows.  We begin with a preliminary section
on periodicity of states of $\F$.
In Section~\ref{section:representations} we construct and analyze
a class of representations of $\On$.
Roughly speaking, this class is indexed by pairs
$(\rho,\theta)$ consisting of a periodic pure state $\rho$ of $\F$
and a representation $\theta$ of $C(\TT)$.
The quasi-orbit of $\rho$ and the unitary equivalence class of
$\theta$ determine the unitary equivalence class of the representation
up to a gauge automorphism.  This ambiguity can be removed with
the addition of a third parameter, called a {\em linking vector\/},
which is related to the periodicity of $\rho$ and is determined 
up to a scalar multiple of modulus one.

In Section~\ref{section:extensions} we study the state extensions
of a periodic pure state $\rho$ of $\F$.
Propositions~\ref{prop:characterize reps} and~\ref{prop:decompose gns}
form the technical core of the paper,
and show that the representations constructed in Section~\ref{section:representations}
include the GNS representation associated with any state which extends $\rho$.
Our main result, Theorem~\ref{theorem:gns}, parameterizes the extensions
of $\rho$ to states of $\On$ by the Borel probability measures on the circle.
In this parameterization the equivalence class of the measure is a complete
invariant for unitary equivalence of state extensions.  We also
compare states which extend different pure states of $\F$.
The invariant we use for this is the set of quasi-equivalence classes
of the shifted states, called the {\em quasi-orbit\/} of $\rho$;
see Definition~\ref{definition:quasiorbit}.
In Corollary~\ref{cor:gns} we answer 
to the affirmative a conjecture made in the final remark of \cite{fowler},
to the effect that a periodic pure state of $\F$ has 
precisely a circle of pure extensions on which the gauge group acts
transitively and $p$-to-$1$, $p$ being the period of the state.
Aperiodic pure states, in contrast, have unique state extensions
which are necessarily pure and fixed by the gauge action
\cite[Theorem 4.3]{laca2}.

In Section~\ref{section:endomorphisms}
we use our representations to construct endomorphisms
of $\bh$ via \eqref{eq:adpi}.  Our main classification result
is Theorem~\ref{theorem:Ad}, where we obtain complete conjugacy invariants
for these endomorphisms based on the parameters $\rho$ and $\theta$.  
In Corollary~\ref{cor:param/classify-endos} we apply
this theorem to classify the endomorphisms 
which arise from extending periodic pure
states of $\F$ to $\On$, as described above.  
The second main result of the section is Theorem~\ref{theorem:tail and fpa},
where we characterize the endomorphisms that arise from state extensions
in terms of their tail and fixed-point algebras.

Although the pure extensions of a pure state $\rho$ are
mutually disjoint, the (ergodic) endomorphisms they produce
are all conjugate.
In Corollary~\ref{cor:classify ergodic} we classify these endomorphisms
using the action of quasi-free automorphisms on the quasi-orbit of $\rho$,
and in Corollary~\ref{cor:characterize ergodic} we characterize them
as those ergodic endomorphisms whose tail algebra has a minimal projection.

In Section~\ref{section:examples} we examine several examples arising from
pure product states of $\F$.  In Example~A we show how our
Theorem~\ref{theorem:gns} generalizes Fowler's result
on pure product states \cite[Theorem~3.1]{fowler}.
In Example~B we consider product states which are constructed from
periodic sequences of unit vectors in $n$\ndash dimensional Hilbert space.
We show that the ergodic endomorphisms which correspond to 
such periodic sequences are completely classified up to
conjugacy by a geometric invariant used in their construction.
This generalizes earlier conjugacy results for shifts from
\cite{powers,sta-shift,laca1,bjp} and for the ergodic endomorphisms
constructed in \cite{laca2}.
Finally, in Example~C we apply our techniques to the problem
of extending the trace on the Choi algebra to $\Oo_2$.

\section{Preliminaries.}\label{prelims}

A {\em multi-index\/} is a $k$\ndash tuple
$s = (s_1,\ldots, s_k)$, where $1 \le s_i \le n$ for each $i$,
and $k$ is any nonnegative integer.  We write $\abs s := k$ and set
$v_s := v_{s_1}\cdots v_{s_k}$, with the convention that $v_s := 1$
if $\abs s = 0$.
Then $\On$ is the closed linear span of monomials of the form
$v_sv_t^*$, where $s$ and $t$ are arbitrary multi-indices,
and $\F$ is the closed linear span of such monomials for which $\abs s = \abs t$.
The canonical conditional expectation $\Phi:\On\to\F$ is given by
\[
\Phi(v_s v_t^*)
  = \begin{cases}
  v_s v_t^* & \text{if $\abs s = \abs t$} \\
  0 & \text{otherwise.} \end{cases}
\]

There are two ways to shift an essential state $\rho$ of $\F$:
`backwards' by $\alpha^*$, as defined in the introduction, and `forwards' by
$\beta^*$, as defined by 
\[
\beta^*\rho(x) = \rho(v_1^* x v_1),\qquad x\in\F.
\]
The arbitrary choice of $v_1$ is irrelevant up to unitary equivalence. 
The shift $\beta^*$ is a quasi-inverse
of $\alpha^*$ in the sense that $\alpha^*\beta^*\rho = \rho \qe \beta^*\alpha^*\rho$
for any essential state $\rho$ of $\F$ \cite[Lemma~3.1]{laca2}.
(We use $\qe$ and $\ue$ to denote quasi-equivalence and unitary equivalence, respectively.)

\begin{example}\label{example:product}
It is helpful to see how the shifts $\alpha^*$ and $\beta^*$ act 
on product states.  Suppose $n$ is finite and $\Ee$ is the
$n$\ndash dimensional Hilbert space spanned by the $v_i$'s,
so that $\Kk(\Ee)$ is isomorphic to the algebra $M_n$ of $n\times n$ matrices.
Then $\F$ is isomorphic to the UHF algebra $M_n \otimes M_n \otimes M_n \dots$
via
$v_s v_t^* \mapsto e_{s_1 t_1} \otimes e_{s_2 t_2} \otimes \dots \otimes e_{s_k t_k}$,
where  $s$ and $t$ are multi-indices of the same length $k$,
and $\{e_{ij}\}$ is the obvious system of matrix units in $M_n$ \cite{cun}.

Suppose $\omega_i$ is a state of $M_n$ for each $i$, and let
\[
\omega = \omega_1 \otimes \omega_2\otimes \omega_3\otimes\cdots 
\]
be the corresponding  product state of $\F$.
Let $\omega_{v_1}$ be the pure
state of $M_n$ determined by $\omega_{v_1}(e_{11}) = 1$. Then
\begin{align*}
\alpha^*\omega 
& = \omega_2 \otimes \omega_3 \otimes \omega_4 \otimes \dotsb \\
\intertext{and}
\beta^*\omega 
& = \omega_{v_1} \otimes \omega_1 \otimes \omega_2 \otimes \dotsb. 
\end{align*}
Similar considerations apply to product states of $\Ff_\infty$ \cite[\S 3]{laca2}.
\end{example}

\begin{definition}\label{definition:quasiorbit}
The {\em quasi-orbit\/} of an essential state $\rho$ of $\F$ is the
set of quasi-equivalence classes
of the states $\alpha^{*k}\rho$ and $\beta^{*k}\rho$ for $k \ge 0$.
\end{definition}

Let us describe the quasi-orbit of an essential factor state $\rho$.
The states $\alpha^{*k}\rho$ and $\beta^{*l}\rho$ for $k,l\ge0$
are factor states (of the same type as $\rho$) \cite[Corollary~3.5]{laca2},
so any given pair of these states is either disjoint or quasi-equivalent.
In the latter case, since both $\alpha^*$ and $\beta^*$ respect
quasi-equivalence of factor states \cite[Corollary~3.6]{laca2}
we can apply an appropriate power of one of the shifts to the
quasi-equivalent pair to obtain $\rho\qe\alpha^{*p}\rho$ for some $p$.

\begin{definition}\label{definition:period}
Suppose $\rho$ is an essential factor state of $\F$.
The {\em period\/} of $\rho$ is the smallest positive integer
$p$ for which $\rho$ is quasi-equivalent to $\alpha^{*p}\rho$.
If no such $p$ exists, we say that $\rho$ is {\em aperiodic\/},
or that it has period $p=\infty$.
\end{definition}

The quasi-orbit of an essential factor state $\rho$ with finite
period $p$ is thus
\[
\set{[\rho],[\alpha^*\rho], \dots, [\alpha^{*(p-1)}\rho]},
\]
or alternatively
\[
\set{[\rho],[\beta^*\rho], \dots, [\beta^{*(p-1)}\rho]},
\]
where the brackets denote quasi-equivalence classes.
In particular, the period of an essential factor state is the cardinality
of its quasi-orbit.

\begin{remark}
Although it would be more accurate to refer to a state which is quasi-equivalent to
its $p^{th}$ translate as {\em quasi-periodic\/}, we will adhere to
the prevailing practice and use the term `periodic' in an asymptotic sense. 
Examples of strictly periodic states (i.e. states which are {\em equal\/}
to their translate by some power of $\alpha^*$) will appear in
Section~\ref{section:examples}.
\end{remark}

Quasi-equivalence of essential factor states of $\F$ is an asymptotic property
(by \cite[Theorem~2.7]{pow-th}
for $n< \infty$ and \cite[Proposition~3.6]{laca1} for $n = \infty$),
so two essential factor states $\rho$ and $\omega$ have the same quasi-orbit
if and only if they are {\em shift-equivalent\/} in the sense that
there exists $k$ such that
\[
\norm{\alpha^{*(k+j)}\rho - \alpha^{*j}\omega} \to 0
\quad \text{ as } \quad j \to \infty.
\]
When $\rho$ and $\omega$ are pure this condition simplifies significantly.
Since $\beta^*$ preserves purity \cite[Lemma~4.2]{laca1},
$\rho$ and $\omega$ have the same quasi-orbit if and only if
\[
\rho\ue\beta^{*k}\omega \quad\text{ or } \quad
\omega\ue\beta^{*k}\rho
\qquad\text{for some $k\ge0$.}
\]
It should be noted that $\alpha^*\rho$ need not be pure even if $\rho$ is.
Some examples of this have been given in \cite{bj}.

We close this section by highlighting some relations between the shifts
of a pure essential state $\rho$ with finite period $p$:
\begin{gather}
\beta^{*k}\rho \ue \beta^{*l}\rho
 \quad \iff \quad
 \alpha^{*k}\rho \qe \alpha^{*l}\rho
 \quad{\iff}\quad
\text{$p$ divides $k - l$, and}
 \label{eq:period} \\
\alpha^{*k}\rho \qe \beta^{*l}\rho
\quad \iff \quad
\text{$p$ divides $k + l$.} \label{eq:period2}
\end{gather}

\section{A class of representations of $\On$.}\label{section:representations}

Suppose $\tilde\rho$ is a state of $\On$, $\rho$ is the restriction
of $\tilde\rho$ to $\F$, and $\tilde\sigma:\On\to\Bb(\tilde\Hh)$
is the GNS representation for $\tilde\rho$ with canonical cyclic vector $\xi$. 
From the unit vectors $\tilde\sigma(v_1^k)\xi$, with $k = 0, 1, 2, \ldots$,
we see that the states $\beta^{*k}\rho$ are vector states
in the restriction of $\tilde\sigma$ to $\F$:
\[
\beta^{*k}\rho(x)
= \ip{\tilde\sigma(x)\tilde\sigma(v_1^k)\xi}{\tilde\sigma(v_1^k)\xi},
\qquad x\in\F.
\]
As an immediate result, the GNS representations of these shifted states  
appear as subrepresentations of the restriction of $\tilde\sigma$ to $\F$. 
Because of this simple fact, whenever we are extending states or
representations from $\F$ to $\On$, we are forced to consider the shifted states. 
It is therefore convenient to establish the following notation,
to be used throughout this paper.
 
\begin{notation}\label{notation:pi}
Suppose $\rho$ is a pure state of $\F$ with finite period
$p$; if $n=\infty$ assume that $\rho$ is essential.
For $i = 0$, $1$, \dots, $p-1$, denote by  $\pi^\rho_i:\F\to\Bb(\Hh^\rho_i)$
the GNS representation for $\beta^{*i}\rho$ with canonical cyclic vector $\xi^\rho_i$.
When there is no chance of confusion we will drop the superscript $\rho$.

For notational convenience we define
$\Hh_p := \Hh_0$ and $\pi_p := \pi_0$.
Our convention for $\xi_p$ will be somewhat different:
\end{notation}

\begin{definition} \label{definition:linking}
A {\em linking vector\/} for $\rho$ is a vector $\xi_p \in \Hh_0$ such that 
\[
\beta^{*p}\rho(x) = \ip{\pi_0(x)\xi_{p}}{\xi_{p}},
\qquad x\in\F.
\]
Since $\beta^{*p}\rho \ue \rho$, there is always a
linking vector, and it is determined up to a scalar multiple of
modulus one because $\rho$ is pure.
\end{definition}

In the following proposition we use a pure essential state $\rho$
and a linking vector $\xi_p$ to construct a representation
$\rep{\rho,\xi_p}$ of $\On$; we will see later (Remark~\ref{remark:irreducible})
that this representation is irreducible.

\begin{prop}\label{prop:S}
Suppose $\rho$ is a pure state of $\F$ with finite period $p$;
if $n=\infty$ assume that $\rho$ is essential.
Let $\xi_p\in\Hh_0$ be a linking vector for $\rho$.

\medskip

\textup{(1)}
If $1 \le a \le n$ and $0 \le i \le p-1$, then
there is an isometry $S_{a,i} : \Hh_i \to \Hh_{i+1}$
determined by
\begin{equation}\label{eq:S}
S_{a,i}\pi_i(x)\xi_i = \pi_{i+1}(v_axv_1^*)\xi_{i+1},
\qquad x\in\F.
\end{equation}

\medskip

\textup{(2)} Let $S_a$ be the isometry $\bigoplus_{i=0}^{p-1} S_{a,i}$
on $\bigoplus_{i=0}^{p-1} \Hh_i$.  There is a representation
$\rep{\rho,\xi_p}$ of $\On$, essential if $n=\infty$, such that
\[
\rep{\rho,\xi_p}(v_a) = S_a,\qquad 1 \le a \le n.
\]

\medskip

\textup{(3)} If $k = i+mp$ with $0 \le i \le p-1$ and $m\ge0$,
then $S_1^k\xi_0$ is a unit vector in $\Hh_i$ which implements
$\beta^{*k}\rho$ as a vector state in $\pi_i$.
For $0 \le k \le p$ we have $S_1^k\xi_0 = \xi_k$, and for
$k\ge p+1$ we define $\xi_k := S_1^k\xi_0$.
\end{prop}

\begin{proof} If $x\in\F$, then
\begin{multline*}
\norm{\pi_{i+1}(v_axv_1^*)\xi_{i+1}}^2
= \ip{\pi_{i+1}(v_1x^*v_a^*v_axv_1^*)\xi_{i+1}}{\xi_{i+1}}
 = \beta^{*(i+1)}\rho(v_1x^*xv_1^*) \\
 = \beta^{*i}\rho(x^*x)
 = \ip{\pi_i(x^*x)\xi_i}{\xi_i}
 = \norm{\pi_i(x)\xi_i}^2.
\end{multline*}
Since vectors of the form $\pi_i(x)\xi_i$ are dense in $\Hh_i$, this gives (1).

We will next show that $S_{a,i}S_{a,i}^* = \pi_{i+1}(v_av_a^*)$,
for which we first need to find a formula for $S_{a,i}^*$.
If $x,y\in\F$, then
\begin{align*}
\ip{S_{a,i}^*\pi_{i+1}(x)\xi_{i+1}}{\pi_i(y)\xi_i}
& = \ip{\pi_{i+1}(x)\xi_{i+1}}{S_{a,i}\pi_i(y)\xi_i} \\
& = \ip{\pi_{i+1}(x)\xi_{i+1}}{\pi_{i+1}(v_ayv_1^*)\xi_{i+1}} \\
& = \ip{\pi_{i+1}(v_1y^*v_a^*x)\xi_{i+1}}{\xi_{i+1}} \\
& = \beta^{*(i+1)}\rho(v_1y^*v_a^*x) \\
& = \beta^{*i}\rho(y^*v_a^*xv_1) \\
& = \ip{\pi_i(y^*v_a^*xv_1)\xi_i}{\xi_i} \\
& = \ip{\pi_i(v_a^*xv_1)\xi_i}{\pi_i(y)\xi_i},
\end{align*}
so
\begin{equation}\label{eq:adjoint}
S_{a,i}^*\pi_{i+1}(x)\xi_{i+1} = \pi_i(v_a^*xv_1)\xi_i,
\qquad x\in\F.
\end{equation}
Using the definition of $S_{a,i}$ we have
\begin{align*}
S_{a,i}S_{a,i}^*\pi_{i+1}(x)\xi_{i+1}
& = S_{a,i}\pi_i(v_a^* x v_1)\xi_i \\
& = \pi_{i+1}(v_a v_a^* x v_1 v_1^*)\xi_{i+1} \\
& = \pi_{i+1}(v_a v_a^*)\pi_{i+1}(x)\pi_{i+1}(v_1 v_1^*)\xi_{i+1},
\end{align*}
so to show that $S_{a,i}S_{a,i}^* = \pi_{i+1}(v_a v_a^*)$ we must verify that
\begin{equation}\label{eq:v1v1*}
\pi_{i+1}(v_1v_1^*)\xi_{i+1} = \xi_{i+1},\qquad 0 \le i \le p-1.
\end{equation}
Since $\pi_{i+1}(v_1v_1^*)$ is a projection, this follows from the calculation
\[
\norm{\pi_{i+1}(v_1v_1^*)\xi_{i+1}}^2
= \ip{\pi_{i+1}(v_1v_1^*)\xi_{i+1}}{\xi_{i+1}}
= \beta^{*(i+1)}\rho(v_1v_1^*)
= \beta^{*i}\rho(1)
= 1.
\]

It is now easy to see that the range projections $S_aS_a^*$ sum to the identity
operator, from which the existence of the representation $\rep{\rho,\xi_p}$
follows immediately:  since each $\pi_i$ is essential,
\[
\sum_{a=1}^n S_aS_a^*
= \sum_{a=1}^n \bigoplus_{i=0}^{p-1} S_{a,i}S_{a,i}^*
= \sum_{a=1}^n \bigoplus_{i=0}^{p-1} \pi_{i+1}(v_av_a^*)
= \bigoplus_{i=0}^{p-1} I_{i+1}
= I.
\]
This completes the proof of (2).

To see that $S_1^k\xi_0$ implements $\beta^{*k}\rho$
as a vector state in $\pi_i$, first observe that
\begin{equation}\label{eq:AdS}
S_1^*\pi_{j+1}(x)S_1 = \pi_j(v_1^*xv_1),\qquad x\in\F,\ 0\le j \le p-1;
\end{equation}
this is an easy consequence of \eqref{eq:S} and \eqref{eq:adjoint}.  Thus
\[
\ip{\pi_i(x)S_1^k\xi_0}{S_1^k\xi_0}
= \ip{S_1^{*k}\pi_i(x)S_1^k\xi_0}{\xi_0}
= \ip{\pi_0(v_1^{*k}xv_1^k)\xi_0}{\xi_0}
= \rho(v_1^{*k}xv_1^k)
= \beta^{*k}\rho(x).
\]
By \eqref{eq:v1v1*} we have $S_1\xi_k = \pi_{k+1}(v_1v_1^*)\xi_{k+1} = \xi_{k+1}$
for $0 \le k \le p-1$, so $\xi_k = S_1^k\xi_0$ for $0 \le k \le p$.
\end{proof}

We are now ready to construct the representations of $\On$ that
will be used to classify the state extensions of $\rho$ to $\On$.
Suppose $U_0$, \dots, $U_{p-1}$ are unitary operators on a Hilbert space $\Kk$.
It is immediate from Proposition~\ref{prop:S} that the range projections of the isometries
$\bigoplus_{i=0}^{p-1} S_{a,i}\otimes U_i$ for $1 \le a \le n$
sum to the identity operator on $\bigoplus_{i=0}^{p-1} \Hh_i\otimes\Kk$.
Consequently there is a unique representation
of $\On$, essential if $n=\infty$,
which maps $v_a$ to $\bigoplus_{i=0}^{p-1} S_{a,i}\otimes U_i$.
We will restrict our attention
to $p$\ndash tuples $(U_0,\dots,U_{p-1})$ 
in which every component but the last one is
equal to the identity; 
up to unitary equivalence of the resulting representation
there is no loss of generality in this, as we will see
in Proposition~\ref{prop:characterize reps}.

\begin{notation}\label{notation:rep}
Consider a triple $(\rho, \xi_p, \theta)$ in which
\begin{itemize}
\item $\rho$ is a pure state of $\F$ with finite period $p$, essential if $n = \infty$,
\item $\xi_p\in\Hh_0$ is a linking vector for $\rho$, and
\item  $\theta$ is a representation of $C(\TT)$ on a Hilbert space $\Kk_\theta$.
\end{itemize}
Let $U_\theta$ be the $p$\ndash tuple 
$(I,I, \ldots, \theta(\mathbf z))$ of unitaries on $\Kk_\theta$, 
where $\mathbf z$ is the identity function on $\TT$.
We will denote by $\rep{\rho,\xi_p,\theta}$ the representation
of $\On$ on $\tilde\Kk_\theta:=\bigoplus_{i=0}^{p-1} \Hh_i\otimes \Kk_\theta$
which is determined by
\begin{equation}\label{eq:rep}
\rep{\rho,\xi_p,\theta}(v_a) = \bigoplus_{i=0}^{p-1} S_{a,i}\otimes U_{\theta,i},
\qquad 1 \le a \le n.
\end{equation}
\end{notation}

\begin{prop}\label{prop:properties}
Suppose $\rho$ is a pure state of $\F$ with finite period $p$;
if $n=\infty$ assume that $\rho$ is essential.
Suppose $\xi_p\in\Hh_0$ is a linking vector for $\rho$ and
$\theta$ is a representation of $C(\TT)$ on a Hilbert space $\Kk_\theta$.

\medskip

\textup{(1)}
The restriction of $\rep{\rho,\xi_p,\theta}$ to $\F$ is
$\bigoplus_{i=0}^{p-1} \pi_i\otimes I_\theta$.

\medskip

\textup{(2)}
If $\psi$ is another representation of $C(\TT)$,
then $\rep{\rho,\xi_p,\theta}$ and $\rep{\rho,\xi_p,\psi}$
are unitarily equivalent (resp. disjoint) if and only if
$\theta$ and $\psi$ are unitarily equivalent (resp. disjoint).

\medskip

\textup{(3)}
$\rep{\rho,\xi_p,\theta}$ is irreducible if and only if $\theta$ is irreducible
(i.e., $\dim\Kk_\theta = 1$).

\medskip

\textup{(4)}
If $\eta\in\Kk_\theta$, then $\xi_0 \otimes \eta$ is cyclic for $\rep{\rho,\xi_p,\theta}$
if and only if $\eta$ is cyclic for $\theta$.

\medskip

\textup{(5)} For each $\lambda\in\TT$ let $\tau_\lambda$ be translation
by $\lambda$ on $C(\TT)$; that is, $\tau_\lambda f(z) = f(\lambda^{-1}z)$
for $f\in C(\TT)$ and $z\in\TT$.  If $\mu^p = \lambda \in\TT$, then
\begin{equation}\label{eq:gauge}
\rep{\rho,\lambda\xi_p,\theta}
\,=\, \rep{\rho,\xi_p,\theta\circ\tau_{\overline\lambda}}
\,\ue\, \rep{\rho,\xi_p,\theta}\circ\gamma_\mu.
\end{equation}
\end{prop}

\begin{proof} (1)  For the moment write $\tilde\pi$ for $\rep{\rho,\xi_p,\theta}$
and $\sigma$ for $\bigoplus_{i=0}^{p-1} \pi_i\otimes I$.
Since both $\sigma$ and the restriction of $\tilde\pi$
to $\F$ are unital representations
and $\F = \clsp\set{v_s v_t^*: \abs s = \abs t}$,
it suffices, by induction, to show that $\tilde\pi(y) = \sigma(y)$
implies that $\tilde\pi(v_jyv_k^*) = \sigma(v_jyv_k^*)$ whenever
$y \in \F$ and  $1 \le j,k \le n$.  If $x\in\F$, $0 \le i \le p-1$,
and $\eta\in\Kk_\theta$, then
\begin{align*}
\tilde\pi(v_j y v_k^*)(\pi_{i+1}(x)\xi_{i+1} \otimes\eta)
& = \tilde\pi(v_j)\sigma(y)
    \tilde\pi(v_k)^*(\pi_{i+1}(x)\xi_{i+1} \otimes\eta) \\
& = \tilde\pi(v_j)\sigma(y)(\pi_i(v_k^* x v_1)\xi_i \otimes U_{\theta,i}^*\eta) \\
& = \tilde\pi(v_j)(\pi_i(y v_k^* x v_1)\xi_i \otimes U_{\theta,i}^*\eta) \\
& = \pi_{i+1}(v_j y v_k^* x v_1 v_1^*)\xi_{i+1} \otimes U_{\theta,i}U_{\theta,i}^*\eta \\
& = \sigma(v_j y v_k^*)(\pi_{i+1}(x)\xi_{i+1} \otimes\eta)
\qquad\text{(by \eqref{eq:v1v1*}).}
\end{align*}

(2) Let $\Ii$ be the intertwining space
\[
\Ii := \setspace{T\in\Bb(\tilde\Kk_\theta,\tilde\Kk_\psi):
T\rep{\rho,\xi_p,\theta}(z) = \rep{\rho,\xi_p,\psi}(z)T
\quad\forall\,z\in\On},
\]
and let
\[
\Ii_0 := \biggl\{\,T = \bigoplus_{i=0}^{p-1} I_i \otimes T_0
\in \Bb(\tilde\Kk_\theta,\tilde\Kk_\psi):
T_0\theta(f) = \psi(f)T_0\quad \forall f\in C(\TT)\,\biggr\}.
\]
We claim that $\Ii = \Ii_0$, from which (2) follows immediately.

As a first step  we describe the space $\Jj\supseteq\Ii$ defined by
\[
\Jj := \setspace{T\in\Bb(\tilde\Kk_\theta,\tilde\Kk_\psi):
T\rep{\rho,\xi_p,\theta}(x) = \rep{\rho,\xi_p,\psi}(x)T
\quad\forall\,x\in\F}.
\]
By (1), $\Jj$ is the set of operators which intertwine
$\bigoplus_{i=0}^{p-1} \pi_i\otimes I_\theta$
and
$\bigoplus_{i=0}^{p-1} \pi_i\otimes I_\psi$.
Since $\pi_0$, \dots, $\pi_{p-1}$ are irreducible  and mutually disjoint,
we have
\[
\Jj = \biggl\{\,T = \bigoplus_{i=0}^{p-1} I_i \otimes T_i:
T_i\in\Bb(\Kk_\theta,\Kk_\psi)\,\biggr\}. 
\]

Suppose that $T = \bigoplus_{i=0}^{p-1} I_i \otimes T_i \in \Jj$.
For notational convenience let $T_p:=T_0$.
If $x\in\F$, $0 \le i \le p-1$ and $\eta\in \Kk_\theta$, then
\begin{equation}\label{eq:Ttheta}
\begin{split}
T\rep{\rho,\xi_p,\theta}(v_a)(\pi_i(x)\xi_i\otimes\eta)
& = T(\pi_{i+1}(v_axv_1^*)\xi_{i+1} \otimes U_{\theta,i}\eta) \\
& =   \pi_{i+1}(v_axv_1^*)\xi_{i+1} \otimes T_{i+1}U_{\theta,i}\eta,
\end{split}
\end{equation}
whereas
\begin{equation}\label{eq:psiT}
\begin{split}
\rep{\rho,\xi_p,\psi}(v_a)T(\pi_i(x)\xi_i \otimes\eta)
& = \rep{\rho,\xi_p,\psi}(v_a)(\pi_i(x)\xi_i \otimes T_i\eta) \\
& = \pi_{i+1}(v_axv_1^*)\xi_{i+1} \otimes U_{\psi,i}T_i\eta.
\end{split}
\end{equation}

Now suppose that $T\in\Ii$.  By \eqref{eq:Ttheta} and \eqref{eq:psiT} we have
\begin{equation}\label{eq:TthetapsiT}
T_{i+1}U_{\theta,i} = U_{\psi,i}T_i,\qquad 0 \le i \le p-1.
\end{equation}
Setting $i=0$, $1$, \dots, $p-2$ gives $T_0 = T_1 = \dots = T_{p-1}$,
and setting $i=p-1$ gives $T_0\theta(\mathbf z) = \psi(\mathbf z)T_0$.
Since $\mathbf z$ generates $C(\TT)$ this implies that $T\in\Ii_0$,
and thus $\Ii\subseteq\Ii_0$.

Conversely, suppose $T_0$ intertwines $\theta$ and $\psi$,
so that $T := \bigoplus_{i=0}^{p-1} I_i \otimes T_0 \in \Ii_0$.
By setting $T_i:=T_0$ for $1 \le i \le p$, we see that
\eqref{eq:TthetapsiT} holds.
By \eqref{eq:Ttheta} and \eqref{eq:psiT} it follows that
$T\rep{\rho,\xi_p,\theta}(v_a) = \rep{\rho,\xi_p,\psi}(v_a)T$ for each $a$,
so that $T\in\Ii$.  Thus $\Ii = \Ii_0$ as claimed, completing the proof of (2).

(3) Setting $\psi = \theta$ in the proof of (2) gives 
\begin{equation}\label{eq:commutant}
\rep{\rho,\xi_p,\theta}(\On)'
= \biggl\{\,T = \bigoplus_{i=0}^{p-1} I_i \otimes T_0:
    T_0\in \theta(C(\TT))'\,\biggr\},
\end{equation}
from which (3) is immediate.

(4) Let $\Mm \subseteq\tilde\Kk_\theta$ be the cyclic subspace for
$\rep{\rho,\xi_p,\theta}$ generated by $\xi_0\otimes\eta$.
The orthogonal projection $P$ of $\tilde\Kk_\theta$ onto $\Mm$
commutes with $\rep{\rho,\xi_p,\theta}(\On)$, so by \eqref{eq:commutant}
there is a projection $P_0\in\theta(C(\TT))'$
such that $P = \bigoplus_{i=0}^{p-1} I_i \otimes P_0$.
On the other hand, $\Mm$ is the closed linear span of vectors of the form
$\rep{\rho,\xi_p,\theta}(v_s v_t^*)(\xi_0\otimes\eta)$,
where $s$ and $t$ are multi-indices.  Given such a vector,
express
$\abs s - \abs t = j + mp$ for $j\in\set{0,\dots,p-1}$ and $m\in\ZZ$.
By \eqref{eq:rep},
\[
\rep{\rho,\xi_p,\theta}(v_s v_t^*)(\xi_0\otimes\eta)
\in\Hh_j \otimes \theta(\mathbf z^m)\eta,
\]
from which it follows that the range of $P_0$
is the closure of $\theta(C(\TT))\eta$.
Assertion~(4) now follows easily.

(5) If  $x\in\F$ and $\eta\in\Kk_\theta$, then
\begin{align*}
\rep{\rho,\lambda\xi_p,\theta}(v_a)(\pi_i(x)\xi_i\otimes\eta)
& = \begin{cases}
    \pi_{i+1}(v_axv_1^*)\xi_{i+1}\otimes\eta 
       & \text{if $0 \le i \le p-2$} \\
    \pi_0(v_axv_1^*)(\lambda\xi_p)\otimes\theta(\mathbf z)\eta
       & \text{if $i = p-1$}
    \end{cases} \\
& = \begin{cases}
    \pi_{i+1}(v_axv_1^*)\xi_{i+1}\otimes\eta
       & \text{if $0 \le i \le p-2$} \\
    \pi_0(v_axv_1^*)\xi_p\otimes\theta(\lambda\mathbf z)\eta
       & \text{if $i = p-1$}
    \end{cases} \\
& = \begin{cases}
    \pi_{i+1}(v_axv_1^*)\xi_{i+1}\otimes\eta
       & \text{if $0 \le i \le p-2$} \\
    \pi_0(v_axv_1^*)\xi_p\otimes\theta\circ\tau_{\overline\lambda}(\mathbf z)\eta
       & \text{if $i = p-1$}
    \end{cases} \\
& = \rep{\rho,\xi_p,\theta\circ\tau_{\overline\lambda}}(v_a)(\pi_i(x)\xi_i\otimes\eta),
\end{align*}
giving the first half of \eqref{eq:gauge}.
Let $T$ be the unitary operator $\bigoplus_{i=0}^{p-1} I_i\otimes \mu^i I_\theta$
on $\tilde\Kk_\theta$.  Then
\begin{align*}
\rep{\rho,\xi_p,\theta}\circ\gamma_\mu(v_a)T(\pi_i(x)\xi_i\otimes\eta)
& = \mu\rep{\rho,\xi_p,\theta}(v_a)(\pi_i(x)\xi_i\otimes\mu^i\eta) \\
& = \begin{cases}
    \pi_{i+1}(v_axv_1^*)\xi_{i+1}\otimes\mu^{i+1}\eta
    & \text{if $0 \le i \le p-2$} \\
    \pi_0(v_axv_1^*)\xi_p\otimes\lambda\theta(\mathbf z)\eta
    & \text{if $i = p-1$}
    \end{cases} \\
& = \begin{cases}
    \pi_{i+1}(v_axv_1^*)\xi_{i+1}\otimes\mu^{i+1}\eta
    & \text{if $0 \le i \le p-2$} \\
    \pi_0(v_axv_1^*)\xi_p\otimes\theta\circ\tau_{\overline\lambda}(\mathbf z)\eta
    & \text{if $i = p-1$}
    \end{cases} \\
& = T\rep{\rho,\xi_p,\theta\circ\tau_{\overline\lambda}}(v_a)(\pi_i(x)\xi_i\otimes\eta),
\end{align*}
from which the second half of \eqref{eq:gauge} follows.
\end{proof}

\begin{remark}\label{remark:irreducible}
The representation
$\rep{\rho,\xi_p}$ of Proposition~\ref{prop:S} is irreducible
because it is unitarily equivalent to 
$\rep{\rho,\xi_p,\varepsilon_1}$, where $\varepsilon_1$ is evaluation at $1\in\TT$,
and $\rep{\rho,\xi_p,\varepsilon_1}$ is irreducible by 
Proposition~\ref{prop:properties}(3).
\end{remark}

\begin{prop}\label{prop:different states}
Suppose $\rho$ and $\omega$ are periodic pure states of $\F$,
essential if $n=\infty$.  If $\rho$ and $\omega$ have the same quasi-orbit
(and hence the same period $p$),
then there are linking vectors $\xi^\rho_p$ and $\xi^\omega_p$
such that $\rep{\rho,\xi^\rho_p,\theta}$ and $\rep{\omega,\xi^\omega_p,\theta}$
are unitarily equivalent for every representation $\theta$ of $C(\TT)$.
\end{prop}

\begin{proof} Define a relation $\approx$ on the pure essential
states of $\F$ with finite period $p$ as follows: $\rho \approx \omega$
if there are linking vectors $\xi^\rho_p$ and $\xi^\omega_p$
such that $\rep{\rho,\xi^\rho_p,\theta} \ue \rep{\omega,\xi^\omega_p,\theta}$
for every representation $\theta$ of $C(\TT)$.
Then $\approx$ is an equivalence relation.  To show transitivity recall 
that the linking
vector for a pure essential state is unique up to a scalar of modulus one
and observe that if $\rep{\rho,\xi^\rho_p,\theta} \ue \rep{\omega,\xi^\omega_p,\theta}$
and $\mu\in\TT$, then by \eqref{eq:gauge}
\[
\rep{\rho,\mu^p\xi^\rho_p,\theta}
\ue
\rep{\rho,\xi^\rho_p,\theta}\circ\gamma_\mu
\ue
\rep{\omega,\xi^\omega_p,\theta}\circ\gamma_\mu
\ue
\rep{\omega,\mu^p\xi^\omega_p,\theta}.
\]

The proof of the proposition is based on the following two claims.

\noindent {\sc Claim~1:} 
If $\rho \ue \omega$, then $\rho \approx \omega$.

\noindent {\sc Claim~2:} $\rho \approx \beta^{*k}\rho$ for every positive integer $k$.

Given the claims, the proof is easy:  if $\rho$ and $\omega$ have the same
quasi-orbit, then $\rho \ue \beta^{*k}\omega$ for some positive integer $k$,
and by the two claims $\rho \approx \beta^{*k}\omega \approx \omega$.
 
\medskip

\noindent {\em Proof of Claim~1:\/} If $\rho \ue \omega$, then there
is a vector $\zeta\in\Hh^\rho_0$ such that
$\omega(x) = \ip{\pi^\rho_0(x)\zeta}\zeta$ for $x\in\F$.
Let $\xi^\rho_p$ be a linking vector for $\rho$,
and let $S_1$ be the isometry on $\bigoplus_{i=0}^{p-1}
\Hh^\rho_i$ defined in Proposition~\ref{prop:S}.
If $0\le i\le p$ and $x\in\F$, then by \eqref{eq:AdS}
\[
\beta^{*i}\omega(x)
 = \omega(v_1^{*i}xv_1^i)
 = \ip{\pi^\rho_0(v_1^{*i}xv_1^i)\zeta}\zeta
 = \ip{S_1^{*i}\pi^\rho_i(x)S_1^i\zeta}\zeta,
 = \ip{\pi^\rho_i(x)S_1^i\zeta}{S_1^i\zeta},
\]
so $S_1^i\zeta$ implements $\beta^{*i}\omega$ as a vector state in $\pi^\rho_i$.
Hence for each $i\in\set{0,\dots,p-1}$ there is a unique
unitary operator $V_i:\Hh^\omega_i \to \Hh^\rho_i$
which intertwines $\pi^\omega_i$ and $\pi^\rho_i$ and maps $\xi^\omega_i$
to $S_1^i\zeta$.  Define $\xi^\omega_p := V_0^*S_1^p\zeta$; then
$\xi^\omega_p$ is a linking vector for $\omega$.  Let $\theta$ be a representation
of $C(\TT)$, and let $V:\bigoplus_{i=0}^{p-1} \Hh^\omega_i \otimes \Kk_\theta
\to \bigoplus_{i=0}^{p-1} \Hh^\rho_i \otimes \Kk_\theta$
be the unitary operator $\bigoplus_{i=0}^{p-1} V_i \otimes I_\theta$.
Then $V$ intertwines $\rep{\omega,\xi^\omega_p,\theta}$
and $\rep{\rho,\xi^\rho_p,\theta}$.  To see this, suppose
$1 \le a \le n$, $0 \le i \le p-1$, $x\in\F$ and $\eta\in\Kk_\theta$.
Using Proposition~\ref{prop:properties}(1) and the convention $V_p := V_0$,
\begin{align*}
V\rep{\omega,\xi^\omega_p,\theta}(v_a)(\pi^\omega_i(x)\xi^\omega_i\otimes\eta)
& = V_{i+1}\pi^\omega_{i+1}(v_axv_1^*)\xi^\omega_{i+1}\otimes U_{\theta,i}\eta \\
& = \pi^\rho_{i+1}(v_axv_1^*)S_1^{i+1}\zeta\otimes U_{\theta,i}\eta \\
& = \rep{\rho,\xi^\rho_p,\theta}(v_axv_1^*)(S_1^{i+1}\zeta\otimes U_{\theta,i}\eta) \\
& = \rep{\rho,\xi^\rho_p,\theta}(v_ax)(S_1^i\zeta\otimes\eta) \\
& = \rep{\rho,\xi^\rho_p,\theta}(v_a)(\pi^\rho_i(x)S_1^i\zeta\otimes\eta) \\
& = \rep{\rho,\xi^\rho_p,\theta}(v_a)V(\pi^\omega_i(x)\xi^\omega_i\otimes\eta).
\end{align*}
This completes the proof of Claim~1.

\medskip

\noindent {\em Proof of Claim~2:\/} It suffices to prove the claim for $k=1$.
Let $\omega := \beta^*\rho$.
Then $\pi^\omega_i = \pi^\rho_{i+1}$,
$\Hh^\omega_i = \Hh^\rho_{i+1}$ and $\xi^\omega_i = \xi^\rho_{i+1}$
for $0 \le i \le p-2$.

Fix a linking vector $\xi^\rho_p$ for $\rho$.
By Proposition~\ref{prop:S}(3), $\xi^\omega_p := \xi^\rho_{p+1}$
is a linking vector for $\omega$.  Let $\theta$ be a representation of $C(\TT)$.
We claim that $\rep{\rho,\xi^\rho_p,\theta} \ue \rep{\omega,\xi^\omega_p,\theta}$.
To construct the intertwining unitary, let $W_0:\Hh^\rho_0\to\Hh^\omega_{p-1}$
be the unique unitary operator which intertwines
$\pi^\rho_0$ and $\pi^\omega_{p-1}$ and maps $\xi^\rho_p$ to $\xi^\omega_{p-1}$.
Let $W:\bigoplus_{i=0}^{p-1} \Hh^\rho_i\otimes\Kk_\theta
    \to\bigoplus_{i=0}^{p-1} \Hh^\omega_i\otimes\Kk_\theta$
be the unitary operator which is the identity from
$\Hh^\rho_{i+1} \otimes \Kk_\theta$ to $\Hh^\omega_i \otimes \Kk_\theta$
for $0 \le i \le p-2$, and $W_0 \otimes\theta(\mathbf z)^*$ from
$\Hh^\rho_0 \otimes \Kk_\theta$ to $\Hh^\omega_{p-1} \otimes \Kk_\theta$.
If $1 \le i \le p$, $x\in\F$ and $\eta\in\Kk_\theta$, then
\[
W(\pi^\rho_i(x)\xi^\rho_i\otimes\eta)
= \begin{cases}
  \pi^\omega_{i-1}(x)\xi^\omega_{i-1}\otimes\eta
    & \text{if $1 \le i \le p-1$} \\
  \pi^\omega_{p-1}(x)\xi^\omega_{p-1}\otimes\theta(\mathbf z)^*\eta
    & \text{if $i = p$,}
  \end{cases}
\]
so for each $a\in{1,\dots,n}$ we have
\begin{align*}
W\rep{\rho,\xi^\rho_p,\theta}(v_a)(\pi^\rho_i(x)\xi^\rho_i\otimes\eta)
& = \begin{cases}
    W(\pi^\rho_{i+1}(v_axv_1^*)\xi^\rho_{i+1}\otimes\eta)
    & \text{if $1 \le i \le p-2$} \\
    W(\pi^\rho_0(v_axv_1^*)\xi^\rho_p\otimes\theta(\mathbf z)\eta)
    & \text{if $i = p-1$} \\
    W(\pi^\rho_1(v_axv_1^*)\xi^\rho_{p+1}\otimes\eta)
    & \text{if $i = p$}
    \end{cases} \\
& = \pi^\omega_i(v_axv_1^*)\xi^\omega_i\otimes\eta \\
& = \begin{cases}
    \rep{\omega,\xi^\omega_p,\theta}(v_a)(\pi^\omega_{i-1}(x)\xi^\omega_{i-1}\otimes\eta)
    & \text{if $1 \le i \le p-1$} \\
    \rep{\omega,\xi^\omega_p,\theta}(v_a)
       (\pi^\omega_{p-1}(x)\xi^\omega_{p-1}\otimes\theta(\mathbf z)^*\eta)
    & \text{if $i = p$}
    \end{cases} \\
& = \rep{\omega,\xi^\omega_p,\theta}(v_a)W(\pi^\rho_i(x)\xi^\rho_i\otimes\eta).
\end{align*}
Thus $W$ intertwines $\rep{\rho,\xi^\rho_p,\theta}$
and $\rep{\omega,\xi^\omega_p,\theta}$, completing the proof of Claim~2.
\end{proof}

\begin{cor}\label{cor:different states}
Suppose $\rho$ and $\omega$ are pure states of $\F$ with finite periods
$p$ and $q$, respectively; if $n=\infty$ assume also that $\rho$ and $\omega$ are
essential.  Suppose $\theta$ and $\psi$ are representations of
$C(\TT)$.  Then $\rep{\rho,\xi^\rho_p,\theta}$ and $\rep{\omega,\xi^\omega_q,\psi}$
are unitarily equivalent (for some choice of linking vectors
$\xi^\rho_p$ and $\xi^\omega_q$) if and only if

\medskip

\textup{(I)} $\rho$ and $\omega$ have the same quasi-orbit, and

\medskip

\textup{(II)} $\theta$ is unitarily equivalent to $\psi\circ\tau_\lambda$
for some $\lambda\in\TT$.
\end{cor}

\begin{proof}
Suppose $\rep{\rho,\xi^\rho_p,\theta} \ue \rep{\omega,\xi^\omega_q,\psi}$.
By Proposition~\ref{prop:properties}(1), the unitary operator which
implements this equivalence also intertwines the representations
$\bigoplus_{i=0}^{p-1} \pi^\rho_i\otimes I_\theta$
and
$\bigoplus_{j=0}^{q-1} \pi^\omega_j\otimes I_\psi$
of $\F$.
From this it is evident that $\Kk_\theta \cong \Kk_\psi$, $p = q$,
and $\pi^\omega_0 \ue \pi^\rho_k$ for some $k\in\set{0,1,\dots,p-1}$.
Thus $\omega \ue \beta^{*k}\rho$, so $\rho$ and $\omega$
have the same quasi-orbit.

By Proposition~\ref{prop:different states}
and the essential uniqueness of linking vectors,
there are scalars $a,b\in\TT$ such that
$\rep{\rho,a\xi^\rho_p,\phi}
 \ue \rep{\omega,b\xi^\omega_p,\phi}$
for every represention $\phi$ of $C(\TT)$.
By \eqref{eq:gauge} we then have
$\rep{\rho,\xi^\rho_p,\phi\circ\tau_{\overline a}}
\ue
\rep{\omega,\xi^\omega_p,\phi\circ\tau_{\overline b}}$
for each $\phi$, and taking  $\phi = \psi\circ\tau_b$ gives
$\rep{\rho,\xi^\rho_p,\psi\circ\tau_{b\overline a}}
 \ue \rep{\omega,\xi^\omega_p,\psi}$.
Thus
$\rep{\rho,\xi^\rho_p,\theta}
 \ue \rep{\rho,\xi^\rho_p,\psi\circ\tau_{b\overline a}}$,
so by Proposition~\ref{prop:properties}(2) we have
$\theta\ue\psi\circ\tau_{b\overline a}$.

Conversely, suppose $\rho$ and $\omega$ have the same quasi-orbit
and $\theta \ue \psi\circ\tau_\lambda$.  By Proposition~\ref{prop:different states},
there are linking vectors $\xi^\rho_p$ and $\xi^\omega_p$ such that
\[
\rep{\rho,\xi^\rho_p,\psi\circ\tau_\lambda}
\ue \rep{\omega,\xi^\omega_p,\psi\circ\tau_\lambda}.
\]
The first of these representations is unitarily equivalent to
$\rep{\rho,\xi^\rho_p,\theta}$ (Proposition~\ref{prop:properties}(2)), 
and the second to $\rep{\omega,\overline\lambda\xi^\omega_p,\psi}$
(Proposition~\ref{prop:properties}(5)).  Thus
$\rep{\rho,\xi^\rho_p,\theta} \ue \rep{\omega,\overline\lambda\xi^\omega_p,\psi}$.
\end{proof}

\section{Extensions of periodic pure states of $\F$ to $\On$.}
\label{section:extensions}

In this section we use the representations constructed in
Section~\ref{section:representations} to parameterize and classify the state extensions
of periodic pure essential states of $\F$.
Our main result, Theorem~\ref{theorem:gns}, is preceded
by a general technical lemma and two technical propositions.

\begin{lemma}\label{lemma}
Suppose $\pi$ is a representation of a \cstar algebra $A$ on a Hilbert space $\Hh$
and $\{\pi_i: i\in I\}$ is a collection of subrepresentations of $\pi$,
each of which is quasi-equivalent to a given representation $\phi$.
Let $\Hh_i$ be the representation space of $\pi_i$.
If $\bigcup_{i\in I} \Hh_i$ has dense linear span in $\Hh$,
then no subrepresentation of $\pi$ is disjoint from $\phi$.
If in addition $\phi$ is factorial, then  $\pi$ is factorial and
 quasi-equivalent to $\phi$.
\end{lemma}

\begin{proof} Suppose $\psi$ is a subrepresentation of $\pi$,
and let $\xi$ be a nonzero vector in the representation space of $\psi$.
Since $\bigcup_{i\in I} \Hh_i$ has dense linear span in $\Hh$,
there is an $i\in I$ such that $\xi \notin \Hh_i^\perp$.
Express $\xi = \xi_0 + \xi_1 \in \Hh_i \oplus \Hh_i^\perp$,
and let $\omega_\xi$, $\omega_{\xi_0}$ and $\omega_{\xi_1}$
be the corresponding vector functionals.
Then $\omega_\xi = \omega_{\xi_0} + \omega_{\xi_1}$,
so $\omega_{\xi_0} \le \omega_\xi$.
Let $\pi_\xi$ and $\pi_{\xi_0}$ denote the GNS representations
associated with $\omega_\xi$ and $\omega_{\xi_0}$, respectively.
Then $\pi_{\xi_0}$ is unitarily equivalent to a subrepresentation
of $\pi_\xi$, which in turn is unitarily equivalent to a subrepresentation
of $\psi$.  But $\pi_{\xi_0}$ is also unitarily equivalent
to a subrepresentation of $\pi_i$, which is quasi-equivalent to $\phi$.
Thus $\psi$ and $\phi$ are not disjoint.
If $\phi$ is factorial, this means that $\pi\qe\phi$. 
\end{proof}

\begin{prop}\label{prop:characterize reps}
Let $\tilde\sigma$ be a representation of $\On$
on a separable Hilbert space $\tilde\Kk$;
if $n=\infty$ assume that $\tilde\sigma$ is essential.
Suppose there exists a pure state $\rho$ of $\F$ with finite period $p$ such that
the restriction of $\tilde\sigma$ to $\F$ decomposes as a direct sum
$\bigoplus_{i=0}^{p-1} \sigma_i$, where $\sigma_i$ is quasi-equivalent to the GNS
representation $\pi_i$ for $\beta^{*i}\rho$. 
Then there is a linking vector $\xi_p$ for
$\rho$ and a representation $\theta$ of $C(\TT)$ such that
$\tilde\sigma$ is unitarily equivalent to $\rep{\rho,\xi_p,\theta}$.
Consequently, the multiplicity of $\pi_i$ in $\sigma_i$
is independent of $i = 0$, $1$, \dots, $p-1$.
\end{prop}

\begin{proof} Let $\sigma = \bigoplus_{i=0}^{p-1} \sigma_i$ be the restriction of
$\tilde\sigma $ to $\F$. Each $\sigma_i$ is unitarily
equivalent to the representation
$\pi_i\otimes I_i$ of $\F$ on $\Hh_i\otimes\Kk_i$ for some separable Hilbert
space $\Kk_i$, so modulo a unitary equivalence we may assume that
$\tilde\Kk = \bigoplus_{i=0}^{p-1} \Hh_i\otimes\Kk_i$
and $\sigma = \bigoplus_{i=0}^{p-1} \pi_i\otimes I_i$.
Let $\xi_p$ be a linking vector for $\rho$, and adopt the
notation convention $\Kk_p := \Kk_0$.  Of course $\pi_p := \pi_0$
and $\Hh_p := \Hh_0$, as usual.

Fix $i\in\set{0,1,\dots,p-1}$ and $\eta\in\Kk_i$.  We claim that
there is a (necessarily unique) vector $U_i\eta \in \Kk_{i+1}$ such that
\begin{equation}\label{eq:Ui}
\tilde\sigma(v_1)(\xi_i \otimes \eta) = \xi_{i+1} \otimes U_i\eta.
\end{equation}
To begin with, note that for any $x\in\F$,
\begin{align*}
\ip{\sigma(x)\tilde\sigma(v_1)(\xi_i \otimes \eta)}
         {\tilde\sigma(v_1)(\xi_i \otimes \eta)}
& = \ip{\sigma(v_1^* x v_1)(\xi_i \otimes \eta)}{\xi_i \otimes \eta} \\
& = \norm\eta^2 \ip{\pi_i(v_1^* x v_1)\xi_i}{\xi_i} \\
& = \norm\eta^2 \beta^{*i}\rho(v_1^* x v_1) \\
& = \norm\eta^2 \beta^{*(i+1)}\rho(x).
\end{align*}
On the other hand, we can express 
$\tilde\sigma(v_1)(\xi_i \otimes \eta)
 = \sum_{j=1}^p \sum_k \zeta_{j,k} \otimes \delta_{j,k}$,
where $\zeta_{j,k} \in \Hh_j$
and for each $j$ the set $\set{\delta_{j,k}}$
is an orthonormal basis for $\Kk_j$.  Then
\begin{align*}
\ip{\sigma(x)\tilde\sigma(v_1)(\xi_i \otimes \eta)}
         {\tilde\sigma(v_1)(\xi_i \otimes \eta)}
& = \sum_{j=1}^p \Bigl\langle
           \sum_k \pi_j(x)\zeta_{j,k} \otimes \delta_{j,k},
           \sum_l \zeta_{j,l} \otimes \delta_{j,l} \Bigr\rangle \\
& = \sum_{j=1}^p \sum_k \ip{\pi_j(x)\zeta_{j,k}}{\zeta_{j,k}}.
\end{align*}
Since $\beta^{*(i+1)}\rho$ is pure, each of the positive linear functionals
$\omega_{\zeta_{j,k}}(x) := \ip{\pi_j(x)\zeta_{j,k}}{\zeta_{j,k}}$
is a multiple of $\beta^{*(i+1)}\rho$.
However, $\omega_{\zeta_{j,k}}$
is also unitarily equivalent to a multiple of $\beta^{*j}\rho$,
because $\pi_j$ is irreducible.
Since the states $\beta^{*j}\rho$ for $j = 1$, $2$, \dots, $p$ are
mutually disjoint, we thus have $\zeta_{j,k} = 0$ unless $j = i+1$.
Moreover, $\omega_{\zeta_{i+1,k}}$ is a scalar multiple of $\beta^{*(i+1)}\rho$
if and only if $\zeta_{i+1,k}$ is a scalar multiple of $\xi_{i+1}$,
so after simplifying and rearranging, the sum
$\sum_{j=1}^p \sum_k \zeta_{j,k} \otimes \delta_{j,k}$
turns out to be an elementary tensor; specifically,
it belongs to the subspace $\xi_{i+1}\otimes \Kk_{i+1}$.
Thus we can define $U_i: \Kk_i \to \Kk_{i+1}$ by \eqref{eq:Ui},
as claimed.

We next claim that $U_i$ is a unitary operator.
It is evident that $U_i$ is linear, and since 
\begin{equation*}
\ip{U_i\eta}{U_i\zeta} = \ip{\tilde\sigma(v_1)(\xi_i \otimes \eta)}
                            {\tilde\sigma(v_1)(\xi_i \otimes \zeta)} 
 = \ip{\xi_i \otimes \eta}{\xi_i \otimes \zeta}
 = \ip\eta\zeta, \qquad \eta,\zeta\in\Kk_i,
\end{equation*}
$U_i$ is an isometry.
To see that $U_i$ is surjective, suppose $\zeta \in \Kk_{i+1}$.
By \eqref{eq:v1v1*} we have
\[
\xi_{i+1} \otimes \zeta 
= \pi_{i+1}(v_1v_1^*)\xi_{i+1} \otimes \zeta
= \sigma(v_1v_1^*)(\xi_{i+1} \otimes \zeta)
= \tilde\sigma(v_1)\tilde\sigma(v_1^*) (\xi_{i+1} \otimes \zeta).
\]
Now  $\tilde\sigma(v_1^*)(\xi_{i+1} \otimes \zeta)$ can be approximated
by a finite sum of vectors of the form $\pi_j(x)\xi_j\otimes\eta$,
where $0 \le j \le p-1$, $x\in\F$ and $\eta\in\Kk_j$,
and for each such vector
\begin{multline}\label{eq:tildesigma}
\tilde\sigma(v_1)(\pi_j(x)\xi_j\otimes\eta)
 = \tilde\sigma(v_1)\sigma(x)(\xi_j\otimes\eta)
 = \sigma(v_1xv_1^*)\tilde\sigma(v_1)(\xi_j\otimes\eta) \\
 = \sigma(v_1xv_1^*)(\xi_{j+1}\otimes U_j\eta)
 = \pi_{j+1}(v_1xv_1^*)\xi_{j+1}\otimes U_j\eta
 \in \Hh_{j+1}\otimes\ran U_j.
\end{multline}
Thus
\[
\xi_{i+1}\otimes\zeta
 = \tilde\sigma(v_1)\tilde\sigma(v_1^*)(\xi_{i+1} \otimes \zeta)
 \in \bigoplus_{j=0}^{p-1} \Hh_{j+1} \otimes \ran U_j,
\]
which shows that $\zeta\in\ran U_i$.  Thus $U_i$ is surjective.

Define unitary operators $T_0$, \dots, $T_{p-1}$ inductively
by $T_0 := U_{p-1}$ and $T_{i+1} = U_iT_i$ for $0 \le i \le p-2$.
Then $T_{p-1}$ is a unitary operator on $\Kk_{p-1}$,
so $\theta(\mathbf z) = T_{p-1}$ determines a representation
$\theta$ of $C(\TT)$ on $\Kk_{p-1}$.  We claim that
$\tilde\sigma$ is unitarily equivalent to $\rep{\rho,\xi_p,\theta}$.
For this, let
$T:\bigoplus_{i=0}^{p-1} \Hh_i\otimes\Kk_{p-1}
 \to\bigoplus_{i=0}^{p-1} \Hh_i\otimes\Kk_i$
be the unitary $\bigoplus_{i=0}^{p-1} I_i \otimes T_i$.
If $1 \le a \le n$, $0\le i \le p-1$, $x\in\F$ and $\eta\in\Kk_{p-1}$,
then by \eqref{eq:tildesigma}
\begin{align*}
\tilde\sigma(v_a)T(\pi_i(x)\xi_i\otimes\eta)
& = \sigma(v_av_1^*)\tilde\sigma(v_1)(\pi_i(x)\xi_i\otimes T_i\eta) \\
& = \sigma(v_av_1^*)(\pi_{i+1}(v_1xv_1^*)\xi_{i+1}\otimes U_iT_i\eta) \\
& = \begin{cases}
    \pi_{i+1}(v_axv_1^*)\xi_{i+1}\otimes T_{i+1}\eta
       & \text{if $0 \le i \le p-2$} \\
    \pi_0(v_axv_1^*)\xi_p\otimes T_0\theta(\mathbf z)\eta
       & \text{if $i = p-1$}
    \end{cases} \\
& = T\rep{\rho,\xi_p,\theta}(v_a)(\pi_i(x)\xi_i\otimes\eta),
\end{align*}
so $T$ intertwines $\rep{\rho,\xi_p,\theta}$ and $\tilde\sigma$.

The multiplicity of $\pi_i$ in $\sigma_i$ is the dimension of $\Kk_i$.
Since each $U_i:\Kk_i\to\Kk_{i+1}$ is unitary, this multiplicity is constant in $i$.
\end{proof}

\begin{prop}\label{prop:decompose gns}
Suppose $\rho$ is a pure state of $\F$ with finite period $p$;
if $n=\infty$ assume that $\rho$ is essential.
Let $\tilde\sigma$ be the GNS representation of
$\On$ corresponding to a state $\tilde \rho$ extending $\rho$.
Then the restriction of $\tilde \sigma$ to $\F$
decomposes as a direct sum $\bigoplus_{i=0}^{p-1}\sigma_i$, with $\sigma_i$
quasi-equivalent to the GNS representation $\pi_i$ of $\beta^{*i}\rho$.
Furthermore, the decomposition is central and the multiplicity
of $\pi_i$ in $\sigma_i$ is independent of $i = 0, 1, \dots, p-1$.
\end{prop}

\begin{proof}
Let $\tilde\Hh$ be the Hilbert space on which $\tilde\sigma$
represents $\On$, and let $\xi\in\tilde\Hh$ be the canonical cyclic
vector which implements $\tilde\rho$ as a vector state.
Let $\sigma$ denote the restriction of $\tilde\sigma$ to $\F$.
For each $k\in\ZZ$ let 
\[
\Gg_k = \set{z\in \On : \gamma_\lambda(z) = \lambda^k z, \lambda\in\TT}.
\]
Notice that $\Gg_0 = \F$ and, in general, that $\Gg_k$ is the 
$k^{th}$ spectral subspace of $\On$ under the action of the gauge group
$\set{\gamma_\lambda: \lambda\in\TT}$.
Define
\[
\Mm_k := \overline{\set{\tilde\sigma(z)\xi: z\in\Gg_k}}.
\]
Then $\Mm_k$ is invariant under $\sigma(\F)$ and $\tilde\Hh =
\clsp\bigcup_{k\in\ZZ} \Mm_k$. Let $\phi_k$ denote 
the subrepresentation of $\sigma$ obtained
by restricting each of the operators $\sigma(x)$ to $\Mm_k$. 
 We claim that
\begin{equation}\label{pmod}
\phi_k \qe \pi_i, 
\end{equation}
 where $i$ is the unique element of
$\set{0,1,\dots,p-1}$ such that $k-i\in p\ZZ$. 
The proof follows \cite[Lemma~3.5]{laca2}.

Suppose $k \ge 0$.  Since $\Gg_k = \F v_1^k$, the vector
$\tilde\sigma(v_1^k)\xi \in \Mm_k$ is cyclic for $\phi_k$.
Moreover, for $x\in\F$ we have
\[
\ip{\phi_k(x)\tilde\sigma(v_1^k)\xi}{\tilde\sigma(v_1^k)\xi}
 = \ip{\sigma(v_1^{*k} x v_1^k)\xi}\xi
 = \beta^{*k}\rho(x).
\]
By the uniqueness of the GNS representation
and \eqref{eq:period} it follows that $\phi_k \ue \pi_k \ue \pi_i$.

Suppose now that $k < 0$.  Using $s$ to denote a multi-index, for $x\in\F$ we have
\[
\alpha^{*\abs k}\rho(x) = \sum_{\abs s = \abs k} \rho(v_s x v_s^*)
= \sum_{\abs s = \abs k} \ip{\sigma(v_s x v_s^*)\xi}\xi
= \sum_{\abs s = \abs k} \ip{\phi_k(x)\tilde\sigma(v_s^*)\xi}{\tilde\sigma(v_s^*)\xi}.
\]
We claim that
$\{\tilde\sigma(v_s^*)\xi: \abs s = \abs k\}$
is generating for $\phi_k$. 
Since $\{v_r v_t^*: \abs r - \abs t = k\}$
has dense linear span in $\Gg_k$, it suffices to show that
\[
\tilde\sigma(v_r v_t^*)\xi \in \clsp\{\sigma(x)\tilde\sigma(v_s^*)\xi:
x\in\F,\,\abs s = \abs k\}
\]
for each such $r$, $t$.  But this is easy: simply write
$t = st'$, where $\abs s = \abs k$, so that
$\tilde\sigma(v_r v_t^*)\xi = \sigma(v_r v_{t'}^*)\tilde\sigma(v_s^*)\xi$.

Since $\{\tilde\sigma(v_s^*)\xi: \abs s = \abs k\}$ is generating for
$\phi_k$, \cite[Lemma~3.2]{laca2} gives that $\phi_k$ is quasi-equivalent
to the GNS representation for $\alpha^{*\abs k}\rho$.
By \eqref{eq:period2},
this implies that $\phi_k \qe \pi_i$, finishing the proof of \eqref{pmod}.

For $i = 0$, $1$, \dots, $p-1$, let $\Ss_i = \clsp\bigcup_{b\in\ZZ} \Mm_{i + bp}$.
Each $\Ss_i$ is invariant under $\sigma(\F)$, and by Lemma~\ref{lemma}
the corresponding subrepresentation $\sigma_i$ of $\sigma$ is quasi-equivalent to $\pi_i$.
The proof  will be complete once we show that
$\Ss_i \perp \Ss_j$ if $i \ne j$, and hence that
$\sigma = \bigoplus_{i=0}^{p-1} \sigma_i$.
For this, suppose $w\in\Gg_k$ and $z\in\Gg_l$,
where $k - l \notin p\ZZ$; we will show that
$\tilde\sigma(w)\xi \perp \tilde\sigma(z)\xi$.
Without loss of generality assume that $k \ge l$.
Let $\zeta = \tilde\sigma(z^*w)\xi \in\Mm_{k - l}$, and write
$\zeta = \zeta_0 \oplus \zeta_1 \in \Mm_0 \oplus \Mm_0^\perp$.
If $\zeta\ne 0$, then the vector functional $\omega_\zeta$
is unitarily equivalent to (a nonzero multiple of) $\beta^{*(k-l)}\rho$.
Since $\beta^{*(k-l)}\rho$ is pure and
$\omega_\zeta  = \omega_{\zeta_0} + \omega_{\zeta_1}$,
either $\omega_{\zeta_0} \ue \beta^{*(k-l)}\rho$
or $\zeta_0 = 0$.  Since $\omega_{\zeta_0} \ue \rho$
and $p$ does not divide $k-l$, by \eqref{eq:period}
we thus have $\zeta_0 = 0$;
that is, $\zeta \perp \Mm_0$.
In particular,
\[
\ip{\tilde\sigma(w)\xi}{\tilde\sigma(z)\xi}
 = \ip\zeta\xi = \ip{\tilde\sigma(z^*w)\xi}\xi
 = 0.
\]

The decomposition $\sigma = \bigoplus_{i=0}^{p-1} \sigma_i$
is central because the $\sigma_i$ are mutually disjoint;
the multiplicity of $\sigma_i$ is constant in $i$ by
Proposition~\ref{prop:characterize reps}.
\end{proof}

\begin{notation}\label{notation:measures}
Let $P(\TT)$ be the space of Borel probability measures on the circle $\TT$.
For each $\mu\in P(\TT)$, let $M_\mu$ be the representation of
$C(\TT)$ on $L^2(\TT,\mu)$ by multiplication operators.
Let $\mathbf 1$ be the function of constant value $1$ on $\TT$.
\end{notation}

\begin{theorem}\label{theorem:gns}
Suppose $\rho$ is a pure state of $\F$ with finite period $p$;
if $n=\infty$ assume that $\rho$ is essential.
For $i\in\set{0,1,\dots,p-1}$ let $\pi_i:\F\to\Bb(\Hh_i)$ be the GNS
representation for $\beta^{*i}\rho$ with canonical cyclic vector $\xi_i$.
Let $\xi_p$ be a linking vector for $\rho$ (as in Definition~\ref{definition:linking}),
and for $k = i+mp\ge p+1$ let $\xi_k$ be the corresponding vector in
$\Hh_i$ which implements $\beta^{*k}\rho$ as a vector state in $\pi_i$
(as in Proposition~\ref{prop:S}(3)).

\medskip

\textup{(1)} For each $\mu\in P(\TT)$ there is a unique state $\tilde\rho[\mu,\xi_p]$
of $\On$ such that
\begin{equation}\label{eq:tilderho}
\tilde\rho[\mu,\xi_p](v_s v_t^*)
= \begin{cases} \ip{\pi_0(v_s v_t^* v_1^{*k})
\xi_k}{\xi_0}\displaystyle\int_{\TT} z^{k/p}\,d\mu(z)
      & \text{if $p$ divides $k$} \\
    0 & \text{otherwise,}
  \end{cases} 
\end{equation}
where $s$ and $t$ are multi-indices with $k: =\abs s - \abs t \ge 0$.
The state $\tilde\rho[\mu,\xi_p]$ extends $\rho$.

\medskip

\textup{(2)} If $\tilde\rho$ is a state of $\On$ which
extends $\rho$, then $\tilde\rho = \tilde\rho[\mu,\xi_p]$ for some $\mu\in P(\TT)$.

\medskip

\textup{(3)} With the linking vector $\xi_p$ fixed,
the map $\mu\mapsto\tilde\rho[\mu,\xi_p]$
is an affine isomorphism of $P(\TT)$ onto the states of $\On$ which extend $\rho$,
and $\tilde\rho[\mu,\xi_p]$ is pure if and only if $\mu$ is a unit point mass.

\medskip

\textup{(4)}
$\tilde\rho[\mu,\xi_p]$ and $\tilde\rho[\nu,\xi_p]$
are unitarily equivalent (resp. disjoint) if and only if the measures $\mu$ and $\nu$
are equivalent (resp. disjoint).

\medskip

\textup{(5)} Suppose $\omega$ is another pure essential state of $\F$.

\qquad\textup{(a)} If $\rho$ and $\omega$ have the same quasi-orbit (and hence the same period $p$),
then there are linking vectors $\xi^\rho_p$ and $\xi^\omega_p$ such that
$\tilde\rho[\mu,\xi^\rho_p]$ and $\tilde\omega[\mu,\xi^\omega_p]$ are unitarily
equivalent for every $\mu\in P(\TT)$.

\qquad\textup{(b)} If $\tilde\rho[\mu,\xi^\rho_p]$ and $\tilde\omega[\nu,\xi^\omega_q]$
are unitarily equivalent, then $\rho$ and $\omega$ have the same quasi-orbit
and $\mu$ is equivalent to a translation of $\nu$.
\end{theorem}

\begin{proof} 
(1) Suppose $s$ and $t$ are multi-indices such that $k:=\abs s - \abs t \ge 0$.
Since elements of the form $v_s v_t^*$ and their adjoints have dense linear span in
$\On$, there is at most one state $\tilde\rho[\mu,\xi_p]$ which satisfies
\eqref{eq:tilderho}. For existence, express $k = i + mp$ with $0 \le i \le p-1$,
and let $\tilde\pi$  be the representation $\rep{\rho,\xi_p,M_\mu}$ of $\On$ on
$\bigoplus_{i=0}^{p-1} \Hh_i\otimes L^2(\TT,\mu)$ defined in \eqref{eq:rep}.  Then
\begin{align*}
\ip{\tilde\pi(v_s v_t^*)(\xi_0\otimes\mathbf 1)}
   {\xi_0\otimes\mathbf 1}
& = \ip{\tilde\pi(v_s v_t^* v_1^{*k}) \tilde\pi( v_1^k)(\xi_0\otimes\mathbf 1)}
        {\xi_0\otimes\mathbf 1} \\
& = \ip{\pi_i(v_s v_t^* v_1^{*k})\xi_k\otimes\mathbf z^m}
        {\xi_0\otimes\mathbf 1} \\
& = \ip{\pi_i(v_s v_t^* v_1^{*k})\xi_k}{\xi_0}\ip{\mathbf z^m}{\mathbf 1} \\
& = \begin{cases}
    \ip{\pi_0(v_s v_t^* v_1^{*k})\xi_k}{\xi_0}\displaystyle\int_{\TT} z^{k/p}\,d\mu(z)
    & \text{if $p$ divides $k$} \\
    0 & \text{otherwise,}
    \end{cases}
\end{align*}
so $\tilde\rho[\mu,\xi_p]$ is the vector state in $\rep{\rho,\xi_p,M_\mu}$
implemented by $\xi_0\otimes\mathbf 1$.  Setting $\abs{s} - \abs{t} = 0 $
shows that $\tilde\rho[\mu,\xi_p]$ extends $\rho$.

(2) Suppose $\tilde\rho$ is a state of $\On$ which extends $\rho$.
By Propositions~\ref{prop:decompose gns} and \ref{prop:characterize reps},
there is a linking vector $\zeta$ for $\rho$ and a representation
$\psi$ of $C(\TT)$ such that the GNS representation for $\tilde\rho$
is unitarily equivalent to $\rep{\rho,\zeta,\psi}$.  Let $\lambda\in\TT$
be such that $\zeta = \lambda\xi_p$,
and let $\theta := \psi\circ\tau_{\overline\lambda}$;
by \eqref{eq:gauge}, $\rep{\rho,\zeta,\psi} \ue \rep{\rho,\xi_p,\theta}$.
Hence there is a unit vector $\xi\in\bigoplus_{i=0}^{p-1} \Hh_i\otimes\Kk_\theta$
which is cyclic for $\rep{\rho,\xi_p,\theta}$ and
which implements $\tilde\rho$ as a vector state.
By the argument used in Proposition~\ref{prop:characterize reps}
to derive \eqref{eq:Ui}, there is a vector $\eta\in\Kk_\theta$
such that $\xi = \xi_0\otimes\eta$.  Thus
\begin{equation}\label{eq:sofar}
\tilde\rho(x) = \ip{\rep{\rho,\xi_p,\theta}(x)(\xi_0\otimes\eta)}{\xi_0\otimes\eta},
\qquad x\in\On.
\end{equation}

Since $\xi_0 \otimes \eta$ is cyclic for $\rep{\rho,\xi_p,\theta}$,
by Proposition~\ref{prop:properties}(4) the vector $\eta$ is cyclic for $\theta$.
It follows that if we define $\mu\in P(\TT)$ by
\[
\int_{\TT} f\,d\mu = \ip{\theta(f)\eta}\eta, \qquad f\in C(\TT),
\]
then $T_0f = \theta(f)\eta$ for $f\in C(\TT)$ determines a unitary operator
$T_0:L^2(\TT,\mu)\to\Kk_\theta$.
Let $T:\bigoplus_{i=0}^{p-1} \Hh_i \otimes L^2(\TT,\mu)
       \to \bigoplus_{i=0}^{p-1} \Hh_i \otimes \Kk_\theta$
be the unitary operator
$\bigoplus_{i=0}^{p-1} I \otimes T_0$.  Routine calculations
show that $T$ intertwines $\rep{\rho,\xi_p,M_\mu}$ and $\rep{\rho,\xi_p,\theta}$
and maps $\xi_0\otimes\mathbf 1$ to $\xi_0 \otimes\eta$.
It now follows immediately from \eqref{eq:sofar} that $\xi_0\otimes\mathbf 1$
implements $\tilde\rho$ as a vector state in $\rep{\rho,\xi_p,M_\mu}$,
so by the proof of (1) we have $\tilde\rho = \tilde\rho[\mu,\xi_p]$.

(3)  From \eqref{eq:tilderho} and part (2)
it is clear that $\mu\mapsto\tilde\rho[\mu,\xi_p]$ is affine and surjective.
To see that it is injective, suppose $\mu,\nu\in P(\TT)$ and $\mu\ne\nu$.
Then there is a positive integer $m$ such that $\int z^m\,d\mu(z) \ne \int z^m\,d\nu(z)$.
Let $k := mp$.
Since $\pi_0$ is irreducible and $\F = \clsp\setspace{v_a v_b^*: \abs a = \abs b}$,
there are multi-indices $a$ and $b$ of equal length such that
$\ip{\pi_0(v_av_b^*)\xi_k}{\xi_0} \ne 0$.  Using $v_j^*v_i = \delta_{ij}1$,
the element $v_av_b^*v_1^k$ can be written in the form $v_sv_t^*$.
Since $\xi_k = \pi_0(v_1^kv_1^{*k})\xi_k$ we then have
\[
\ip{\pi_0(v_sv_t^*v_1^{*k})\xi_k}{\xi_0}
= \ip{\pi_0(v_av_b^*)\pi_0(v_1^kv_1^{*k})\xi_k}{\xi_0}
= \ip{\pi_0(v_av_b^*)\xi_k}{\xi_0}
\ne 0.
\]
By \eqref{eq:tilderho} it follows that
$\tilde\rho[\mu,\xi_p](v_s v_t^*) \ne \tilde\rho[\nu,\xi_p](v_s v_t^*)$,
completing the proof of injectivity.
Since $\mu\mapsto\tilde\rho[\mu,\xi_p]$ is an affine isomorphism
it preserves extreme points; hence point masses correspond to pure states.

(4) Since $\mathbf 1$ is cyclic for $M_\mu$, by Proposition~\ref{prop:properties}(4)
the vector $\xi_0\otimes\mathbf 1$ is cyclic for $\rep{\rho,\xi_p,M_\mu}$.
Thus $\rep{\rho,\xi_p,M_\mu}$ is unitarily equivalent to the
GNS representation for $\tilde\rho[\mu,\xi_p]$.
Since $M_\mu$ and $M_\nu$ are unitarily equivalent (resp. disjoint)
if and only if $\mu$ and $\nu$ are equivalent (resp. disjoint) measures,
(4) follows immediately from Proposition~\ref{prop:properties}(2).

Finally, since $\mu$ is equivalent to a translate of $\nu$ if and only if
$M_\mu \ue M_\nu\circ\tau_\lambda$ for some $\lambda\in\TT$,
(5) is an immediate consequence of Proposition~\ref{prop:different states}
and Corollary~\ref{cor:different states}.
\end{proof}

\begin{cor}\label{cor:gns}
Suppose $\rho$ is a pure state of $\F$ with finite period $p$; if $n=\infty$
assume that $\rho$ is essential.
The gauge group acts $p$-to-$1$ and transitively on the 
pure extensions of $\rho$ to $\On$,
and distinct pure extensions are disjoint.
\end{cor}

\begin{proof}
Fix a linking vector $\xi_p$ for $\rho$.  By Theorem~\ref{theorem:gns}(3),
the pure extensions of $\rho$ are $\setspace{\tilde\rho[\mu_c,\xi_p]: c\in\TT}$,
where $\mu_c$ denotes the unit point mass at $c$.
Since no two different point masses are equivalent, it follows
from Theorem~\ref{theorem:gns}(4) that no two different pure extensions are unitarily
equivalent; that is, distinct pure extensions are disjoint.

If $s$ and $t$ are multi-indices with $k: =\abs s - \abs t \ge 0$,
then by \eqref{eq:tilderho}
\[
\tilde\rho[\mu_c,\xi_p](v_s v_t^*)
= \begin{cases} c^{k/p}\ip{\pi_0(v_s v_t^* v_1^{*k})\xi_k}{\xi_0}
    & \text{if $p$ divides $k$} \\
  0 & \text{otherwise.}
\end{cases}
\]
On the other hand, if $\lambda\in\TT$, then
\[
\tilde\rho[\mu_c,\xi_p]\circ\gamma_\lambda(v_s v_t^*)
= \begin{cases} (\lambda^pc)^{k/p}\ip{\pi_0(v_s v_t^* v_1^{*k})\xi_k}{\xi_0}
    & \text{if $p$ divides $k$} \\
  0 & \text{otherwise,}
\end{cases}
\]
so $\tilde\rho[\mu_c,\xi_p]\circ\gamma_\lambda
  = \tilde\rho[\mu_{\lambda^pc},\xi_p]$.
Thus the gauge action is transitive on the pure extensions of $\rho$,
and for any pure extension $\tilde\rho$ we have
$\tilde\rho\circ\gamma_\lambda = \tilde\rho$ if and only if $\lambda^p = 1$.
\end{proof}

\section{Endomorphisms of $\bh$.}\label{section:endomorphisms}
We are now ready to construct and classify endomorphisms of $\bh$
using the representations from Section~\ref{section:representations}.
We will use $\conj$ to denote conjugacy of endomorphisms.

Recall that a representation $\phi:\On\to\bh$ gives rise
to an endomorphism of $\bh$ via
\[
\Ad\phi(A) = \sum_{k=1}^n \phi(v_k)A\phi(v_k)^*,\qquad A\in\bh.
\]
Recall also that the gauge action $\gamma:\TT\to\Aut(\On)$ extends to an action
of the unitary group $\Uu(\Ee)$ by {\em quasi-free\/} automorphisms, determined by
$\gamma_W(v_a) = Wv_a$.
Modifying $\phi$ by a quasi-free automorphism
does not change $\Ad\phi$, and modifying it by a unitary equivalence
only changes $\Ad\phi$ to a conjugate endomorphism.
This is indeed all the collapsing there is in the map $\phi \mapsto \Ad\phi$:
the endomorphisms $\Ad\phi_1$ and $\Ad\phi_2$ are conjugate if and only if 
$\phi_2 \ue \phi_1 \circ \gamma_W$ for some $W\in \Uu(\Ee)$
\cite[Proposition~2.4]{laca1}.

Suppose $\rho$ is a periodic pure essential state of $\F$
and $\theta$ is a representation of $C(\TT)$.
For each choice of linking vector $\xi_p$ for $\rho$
we can form the representation $\rep{\rho,\xi_p,\theta}$
as in \eqref{eq:rep}.
By Proposition~\ref{prop:properties}(5), two
representations of the form $\rep{\rho,*,\theta}$ differ by at most a
gauge automorphism and a unitary equivalence,
so the conjugacy class of $\Ad\rep{\rho,\xi_p,\theta}$
does not depend on the choice of $\xi_p$.
We will denote this conjugacy class, or a representative thereof, by 
\[
\alpha_{\rho,\theta} := \Ad\rep{\rho,\xi_p,\theta}.
\]
Since we only look at endomorphisms modulo conjugacy,
this slight abuse of notation will not cause problems.

Two endomorphisms coming from different states 
of $\F$ and representations of $C(\TT)$
can be conjugate, and the following theorem determines exactly when this happens. 

\begin{theorem}\label{theorem:Ad}
Suppose $\rho$ and $\omega$ are periodic pure states of $\F$, essential if $n=\infty$,
and let
$\theta$ and $\psi$ be representations of $C(\TT)$. Then $\alpha_{\rho,\theta}$ and
$\alpha_{\omega,\psi}$ are conjugate if and only if

\medskip

\textup{(I)} there is a unitary operator $W$ on $\Ee$ such that
$\rho\circ\gamma_W$ and $\omega$ have the same quasi-orbit, and

\medskip

\textup{(II)} $\theta$ is unitarily equivalent to $\psi\circ\tau_\lambda$
for some $\lambda\in\TT$.
\end{theorem}

\begin{proof} Let $\xi^\rho_p$ and $\xi^\omega_q$ be linking vectors for $\rho$
and $\omega$, respectively.  The endomorphisms
$\alpha_{\rho,\theta}$ and $\alpha_{\omega,\psi}$ are conjugate if and only if
there is a unitary operator $W$ on $\Ee$ such that
$\rep{\rho,\xi^\rho_p,\theta}\circ\gamma_W$ and $\rep{\omega,\xi^\omega_q,\psi}$
are unitarily equivalent.
The proof will be by direct application of Corollary~\ref{cor:different states} 
once $\rep{\rho,\xi^\rho_p,\theta}\circ\gamma_W$ has been changed
to an appropriate form.  It suffices to prove the following:

\medskip
\noindent {\sc Claim:} 
For every unitary $W$ on $\Ee$ there is a linking vector
$\xi^{\rho\circ\gamma_W}_p$ for $\rho\circ\gamma_W$ such that
\[
\rep{\rho,\xi^\rho_p,\theta}\circ\gamma_W
\ue
\rep{\rho\circ\gamma_W,\xi^{\rho\circ\gamma_W}_p,\theta}.
\]

\noindent {\em Proof of Claim:\/}
Fix a unitary operator $W$ on $\Ee$,
let $\xi^\rho_p$ be a linking vector for $\rho$,
and let $S_1$ be the isometry on $\bigoplus_{i=0}^{p-1}
\Hh^\rho_i$ defined in Proposition~\ref{prop:S}.
If $0 \le i \le p$ and $x\in\F$, then by \eqref{eq:AdS}
and Proposition~\ref{prop:S}(3)
\begin{align*}
(\beta^{*i}(\rho\circ\gamma_W))(x)
& = \rho\circ\gamma_W(v_1^{*i}xv_1^i) \\
& = \ip{\pi^\rho_0\circ\gamma_W(v_1^{*i}xv_1^i)\xi^\rho_0}{\xi^\rho_0} \\
& = \ip{S_1^{*i}\pi^\rho_i(v_1^i\gamma_W(v_1^{*i}xv_1^i)v_1^{*i})S_1^i\xi^\rho_0}{\xi^\rho_0} \\
& = \ip{\pi^\rho_i(v_1^i\gamma_W(v_1^{*i}xv_1^i)v_1^{*i})\xi^\rho_i}{\xi^\rho_i} \\
& = \ip{\pi^\rho_i\circ\gamma_W(x)\pi^\rho_i(\gamma_W(v_1^i)v_1^{*i})\xi^\rho_i}
       {\pi^\rho_i(\gamma_W(v_1^i)v_1^{*i})\xi^\rho_i},
\end{align*}
so $\pi^\rho_i(\gamma_W(v_1^i)v_1^{*i})\xi^\rho_i$ implements
$\beta^{*i}(\rho\circ\gamma_W)$ as a vector state in $\pi^\rho_i\circ\gamma_W$.
Hence for each $i\in\set{0,\dots,p-1}$ there is a unique
unitary operator $V_i:\Hh^{\rho\circ\gamma_W}_i \to \Hh^\rho_i$
which intertwines $\pi^{\rho\circ\gamma_W}_i$ and
$\pi^\rho_i\circ\gamma_W$ and satisfies
\[
V_i\xi^{\rho\circ\gamma_W}_i = \pi^\rho_i(\gamma_W(v_1^i)v_1^{*i})\xi^\rho_i.
\]
Define $\xi^{\rho\circ\gamma_W}_p := V_0^*\pi^\rho_0(\gamma_W(v_1^p)v_1^{*p})\xi^\rho_p$;
then $\xi^{\rho\circ\gamma_W}_p$ is a linking vector for $\rho\circ\gamma_W$.
Let $\theta$ be a representation of $C(\TT)$, and let
$V:\bigoplus_{i=0}^{p-1} \Hh^{\rho\circ\gamma_W}_i \otimes \Kk_\theta
\to \bigoplus_{i=0}^{p-1} \Hh^\rho_i \otimes \Kk_\theta$
be the unitary operator $\bigoplus_{i=0}^{p-1} V_i \otimes I_\theta$.
Then $V$ intertwines $\rep{\rho\circ\gamma_W,\xi^{\rho\circ\gamma_W}_p,\theta}$
and $\rep{\rho,\xi^\rho_p,\theta}\circ\gamma_W$.  To see this, suppose
$1 \le a \le n$, $0 \le i \le p-1$, $x\in\F$ and $\eta\in\Kk_\theta$.
Using Proposition~\ref{prop:properties}(1) and the convention $V_p := V_0$,
\begin{align*}
V\rep{\rho\circ\gamma_W,\xi^{\rho\circ\gamma_W}_p,\theta}(v_a)
& (\pi^{\rho\circ\gamma_W}_i(x)\xi^{\rho\circ\gamma_W}_i\otimes\eta) \\
& = V_{i+1}\pi^{\rho\circ\gamma_W}_{i+1}(v_axv_1^*)\xi^{\rho\circ\gamma_W}_{i+1}
  \otimes U_{\theta,i}\eta \\
& = \pi^\rho_{i+1}\circ\gamma_W(v_axv_1^*)
    \pi^\rho_{i+1}(\gamma_W(v_1^{i+1})v_1^{*(i+1)})\xi^\rho_{i+1}
  \otimes U_{\theta,i}\eta \\
& = \pi^\rho_{i+1}(\gamma_W(v_axv_1^i)v_1^{*(i+1)})\xi^\rho_{i+1}
  \otimes U_{\theta,i}\eta \\
& = \rep{\rho,\xi^\rho_p,\theta}(\gamma_W(v_axv_1^i)v_1^{*(i+1)})
    \rep{\rho,\xi^\rho_p,\theta}(v_1)(\xi^\rho_i\otimes\eta) \\
& = \rep{\rho,\xi^\rho_p,\theta}(\gamma_W(v_axv_1^i)v_1^{*i})(\xi^\rho_i\otimes\eta) \\
& = \rep{\rho,\xi^\rho_p,\theta}\circ\gamma_W(v_a)
    (\pi^\rho_i\circ\gamma_W(x)\pi^\rho_i(\gamma_W(v_1^i)v_1^{*i})\xi^\rho_i\otimes\eta) \\
& = \rep{\rho,\xi^\rho_p,\theta}\circ\gamma_W(v_a)V
    (\pi^{\rho\circ\gamma_W}_i(x)\xi^{\rho\circ\gamma_W}_i\otimes\eta).
\end{align*}
This completes the proof of the claim, and hence the theorem.
\end{proof}

When $\theta$ is the representation $M_\mu$ by multiplication
operators on $L^2(\TT,\mu)$ we write simply $\alpha_{\rho,\mu}$
in place of $\alpha_{\rho,M_\mu}$.
As an immediate corollary to Theorem~\ref{theorem:gns}
we can now parameterize and classify the endomorphisms constructed using
the strategy of \cite{laca1}, wherein one starts with a pure essential state
$\rho$ of $\F$, extends $\rho$ to a state $\tilde\rho$ of $\On$, 
and then uses the GNS representation for $\tilde\rho$ to implement
an endomorphism of $\bh$.

\begin{cor}\label{cor:param/classify-endos}
Suppose $\rho$ is a periodic pure state of $\F$, essential if $n=\infty$,
and $\tilde\sigma$ is the GNS representation for
a state $\tilde\rho$ of $\On$ which extends $\rho$.
Then there is a Borel probability measure $\mu$ on the circle $\TT$
such that $\Ad\tilde\sigma$ is conjugate to $\alpha_{\rho,\mu}$.

Let $\omega$ be another periodic pure essential state of $\F$,
and let $\nu\in P(\TT)$.
Then $\alpha_{\rho,\mu}$ and $\alpha_{\omega,\nu}$ are conjugate
if and only if

\medskip

\textup{(I)} there is a unitary operator $W$ on $\Ee$ such that
$\rho\circ\gamma_W$ and $\omega$ have
the same quasi-orbit, and

\medskip

\textup{(II)} $\mu$ is equivalent to a translate of $\nu$.
\end{cor}

\begin{proof}  Fix a linking vector $\xi_p$ for $\rho$.  By
Theorem~\ref{theorem:gns}(2), $\tilde\rho = \tilde\rho[\mu,\xi_p]$
for some $\mu\in P(\TT)$, so that $\tilde\sigma \ue \rep{\rho,\xi_p,M_\mu}$.
Thus $\Ad\tilde\sigma \conj \alpha_{\rho,\mu}$.
 
Since the measure $\mu$ is equivalent to a translate of $\nu$ if and only if
$M_\mu \ue M_\nu\circ\tau_\lambda$ for some $\lambda\in\TT$, the second part follows
directly from Theorem~\ref{theorem:Ad}.
\end{proof}

There are two von Neumann algebras naturally associated with an endomorphism $\alpha$
of $\bh$:  its {\em tail algebra\/}
\[
\tail(\alpha) := \bigcap_{k=1}^\infty \alpha^k(\bh),
\]
and its {\em fixed-point algebra\/}
\[
\fpa(\alpha) := \setspace{A\in\bh:\alpha(A) = A}.
\]
Assuming as we are that $\alpha$ is unital, one can always realize
$\alpha$ as $\Ad\phi$ for some (essential) representation
$\phi$ of $\On$.  By \cite[Proposition~3.1]{laca1},
$\fpa(\alpha)$ is the commutant of $\phi(\On)$
and $\tail(\alpha)$ is the commutant of $\phi(\F)$.
If $\phi$ is the GNS representation of some state $\tilde\rho$ of $\On$,
then the canonical cyclic vector $\xi$ for $\phi$ is separating for
$\fpa(\alpha)$.  In the following theorem we show that when the restriction
of $\tilde\rho$ to $\F$ is pure, much more is true:
$\fpa(\alpha)$ is abelian, and the projection
onto the closure of $\tail(\alpha)'\xi$ is minimal in $\tail(\alpha)$.
Moreover, the latter condition characterizes these endomorphisms.

\begin{theorem}\label{theorem:tail and fpa}
Suppose $\alpha$ is a unital endomorphism of $\bh$ with Powers index
$n$ $(2 \le n \le \infty)$.  Then \textup{(I)} and \textup{(II)}
below are equivalent:

\medskip

\textup{(I)}
$\alpha$ is conjugate to  $\Ad\tilde\sigma$ for $\tilde\sigma$
the GNS representation of a state extending a pure essential state $\rho$ 
of $\F$.

\medskip

\textup{(II)}
$\tail(\alpha)$ has a minimal projection whose
range contains a separating vector for $\fpa(\alpha)$.

\medskip

If $\alpha$ satisfies \textup{(I)} and \textup{(II)},
then the center of $\tail(\alpha)$ is finite-dimensional if and only
if $\rho$ has finite period, in which case $\dim Z(\tail(\alpha))$
is the period of $\rho$. Moreover,

\textup{(1)} If $\rho$ has finite period $p$, then there is a
Borel probability measure $\mu$ on the circle $\TT$ such that
$\tail(\alpha)$ is spatially isomorphic to
\begin{equation}\label{eq:tail}
\biggl\{\,\bigoplus_{i=0}^{p-1} I_i\otimes T_i:
T_i\in\Bb(L^2(\TT,\mu))\,\biggr\}
\subset \Bb\biggl(\bigoplus_{i=0}^{p-1}\Hh_i\otimes L^2(\TT,\mu)\biggr) 
\end{equation}
and $\fpa(\alpha)$ is spatially isomorphic to the abelian algebra
\begin{equation}\label{eq:fpa}
\biggl\{\,\bigoplus_{i=0}^{p-1} I_i \otimes T_0:
T_0\in M_\mu(L^\infty(\TT,\mu))\,\biggr\}
\subset \Bb\biggl(\bigoplus_{i=0}^{p-1}\Hh_i\otimes L^2(\TT,\mu)\biggr),
\end{equation}
where as usual $\Hh_i$ denotes the GNS Hilbert space for $\beta^{*i}\rho$.

\textup{(2)} If $\rho$ is aperiodic, then $\tail(\alpha)$ is isomorphic to
$\ell^\infty(\ZZ)$.
\end{theorem}

\begin{proof} (II) $\Ra$ (I): Suppose $\tail(\alpha)$ has a minimal projection $P$ whose
range contains a vector $\xi$ which is separating for $\fpa(\alpha)$.
Let $\phi$ be a representation of $\On$ such that $\alpha = \Ad\phi$,
and let $\tilde\rho$ be the vector state of $\On$ in $\phi$ implemented by $\xi$.
Since $\xi$ is separating for $\fpa(\alpha) = \phi(\On)'$
it is cyclic for $\phi$, so $\alpha\conj\Ad\tilde\sigma$
with $\tilde\sigma$ the GNS representation for $\tilde\rho$.
The restriction of $\tilde\rho$ to $\F$ is pure because
$P$ is minimal in $\tail(\alpha) = \phi(\F)'$.

(I) $\Ra$ (II): If (I) holds and $\rho$ has finite period $p$,
then by Corollary~\ref{cor:param/classify-endos} there is a Borel probability measure
$\mu$ on $\TT$ such that $\alpha \conj \alpha_{\rho,\mu}$.
Let $\xi_p$ be a linking vector for $\rho$.
The tail and fixed-point algebras of $\alpha$ are then spatially equivalent
to those of $\Ad\rep{\rho,\xi_p,M_\mu}$,
which by the proof of Proposition~\ref{prop:properties}(2) are given by
\eqref{eq:tail} and \eqref{eq:fpa}, respectively; the second of these
requires the extra observation that $M_\mu(C(\TT))' = M_\mu(L^\infty(\TT,\mu))$.
Let $P_0$ be the rank-one projection onto the constant function
$\mathbf 1\in L^2(\TT,\mu)$, let $P_i = 0$ for $i = 1$, \dots, $p-1$,
and let $P = \bigoplus_{i=0}^{p-1} I_i \otimes P_i$; then $P$ is a minimal
projection in the tail algebra.  Since $\mathbf 1$ is cyclic for $M^\mu(C(\TT))$
it is separating for $M^\mu(C(\TT))'$, so any nonzero vector in the range of $P$
is separating for the fixed-point algebra.

If $\rho$ is aperiodic, then $\tilde\rho$ must be the gauge-invariant state
$\rho\circ\Phi$.  Let $\sigma$ be the restriction
of $\tilde\sigma$ to $\F$.  By \cite[Propositions~2.2, 3.4]{laca2},
$\sigma$ decomposes as a direct sum $\bigoplus_{i=-\infty}^\infty \sigma_i$,
where $\sigma_i$ is irreducible and quasi-equivalent to
the GNS representation of $\beta^{*i}\rho$ (resp. $\alpha^{*\abs i}\rho$)
if $i\ge 0$ (resp. $i<0$).  Since these irreducible summands are mutually
disjoint,
\[
\tail(\Ad\tilde\sigma) = \sigma(\F)' \cong \ell^\infty(\ZZ).
\]
Let $P$ be any minimal projection in this algebra.
Any nonzero vector in the range of $P$ is separating for $\fpa(\alpha)$
since $\alpha$ is ergodic.
\end{proof}

\subsection*{Ergodic endomorphisms}
By \cite[Proposition~3.1]{laca1},
$\Ad\phi$ is ergodic if and only if $\phi$ is irreducible.
Thus the pure extensions of a pure essential state
$\rho$ yield ergodic endomorphisms via their GNS representations. 
Since these pure extensions are in the same gauge orbit,
the corresponding endomorphisms are all conjugate:

\begin{cor}\label{cor:classify ergodic}
Suppose $\rho$ is a pure state of $\F$, essential if $n=\infty$.

\medskip

\textup{(1)}
Let $\tilde\sigma$ be the GNS representation for a pure extension of $\rho$.
Then the ergodic endomorphism $\Ad\tilde\sigma$
depends only on $\rho$ up to conjugacy,
so we denote it by $\alpha_\rho:= \Ad \tilde\sigma$.  

\medskip

\textup{(2)}
If $\omega$ is another pure essential state of $\F$, then
$\alpha_\rho$ and $\alpha_\omega$ are conjugate if
and only if there is a unitary operator $W$ on $\Ee$ such that
$\rho\circ\gamma_W$ and $\omega$ have the same quasi-orbit.
\end{cor}

\begin{proof} First suppose $\rho$ is periodic.  Let $\xi_p$
be a linking vector for $\rho$.
If $\tilde\rho$ is a pure state of $\On$ which extends $\rho$,
then by Theorem~\ref{theorem:gns}(3)
there is a unit point mass $\mu$ on $\TT$ such that
$\tilde\rho = \tilde\rho[\mu,\xi_p]$.
The GNS representation $\tilde\sigma$ for $\tilde\rho$
is thus unitarily equivalent to $\rep{\rho,\xi_p,M_\mu}$,
so that $\Ad\tilde\sigma$ is conjugate to $\alpha_{\rho,\mu}$.
Since condition (II) of Corollary~\ref{cor:param/classify-endos}
is automatic for point masses, all such endomorphisms
$\Ad\tilde\sigma$ are conjugate.
The second assertion also follows from Corollary~\ref{cor:param/classify-endos}.

Suppose now that $\rho$ is aperiodic.  By \cite[Theorem~4.3]{laca2},
the gauge-invariant extension $\rho\circ\Phi$
is the only state of $\On$ which extends $\rho$, and
$\rho\circ\Phi$ is pure, so (1) is trivial.
Part (2) follows from \cite[Theorem~4.2]{laca2}.
\end{proof}

Finally we give an intrinsic characterization of the class of ergodic endomorphisms
arising from pure states of $\F$ in terms of the tail algebra.

\begin{cor}\label{cor:characterize ergodic}

\textup{(1)} Suppose $\rho$ is a pure state of $\F$, essential if $n=\infty$.
Let $\alpha_\rho$ be the ergodic endomorphism associated with
$\rho$ as in Corollary~\ref{cor:classify ergodic}.
Then $\tail(\alpha_\rho)$ is isomorphic to $\CC^p$
if $\rho$ has finite period $p$, and $\ell^\infty(\ZZ)$
if $\rho$ is aperiodic.

\medskip

\textup{(2)} Suppose $\alpha$ is an ergodic endomorphism of $\bh$
whose tail algebra has a minimal projection. Then there is a pure
essential state $\rho$ of $\F$ such that $\alpha$ is conjugate to $\alpha_\rho$.
In particular, if $\alpha$ is a shift then it is conjugate to $\alpha_\rho$ 
for some pure essential quasi-invariant state $\rho$.
\end{cor}

\begin{proof} (1) By definition,
$\alpha_\rho$ satisfies condition (I) of Theorem~\ref{theorem:tail and fpa}.
The result is thus immediate from this theorem for aperiodic $\rho$.
If $\rho$ has finite period $p$ then 
$\tail(\alpha_\rho)$ is given by \eqref{eq:tail} for some point mass $\mu$,
and is hence isomorphic to $\CC^p$.

(2) Let $P$ be a minimal projection in $\tail(\alpha)$.
Every nonzero vector in the range of $P$ is separating
for $\fpa(\alpha) = \CC I$, so by (II)$\Ra$(I) of
Theorem~\ref{theorem:tail and fpa} and 
Corollary~\ref{cor:classify ergodic}(1), $\alpha = \alpha_\rho$
for some pure essential state $\rho$ of $\F$. 
If $\alpha$ is a shift then $\tail(\alpha)$ consists of scalar operators
so the identity is a minimal projection. Thus $\alpha = \alpha_\rho$
for some pure essential state $\rho$,
and by \cite[Theorem 4.5]{laca1} $\rho$ must be quasi-invariant.
\end{proof}

We finish the section by pointing out that, as a consequence of the Corollary,
there is an interesting
restriction on the possible tail algebras of ergodic endomorphisms:

\begin{scholium} 
If the tail algebra of an ergodic endomorphism
has a minimal projection, then it is isomorphic to either $\CC^p$
or $\ell^\infty(\ZZ)$, depending on the period $p$ of the state
arising from a vector in the range of the minimal projection. 
\end{scholium}

\section{Examples}\label{section:examples}
Our main source of examples are the pure {\em product states\/}
$
\omega = \otimes_{i=1}^\infty \omega_i
$
of $\F$, where each $\omega_i$ is a pure state of $\Kk(\Ee)$;
c.f. Example~\ref{example:product}.
For each unit vector $v$ in $\Ee$ let $\omega_v$ be the pure
state of $\Kk(\Ee)$ given by $\omega_v (T) = \langle T v,v\rangle$;
strictly speaking,  $\omega_v$  depends only on the 
one-dimensional subspace $[v] := \mathbb C v$ and not on $v$ itself.
If $f = (f_1, f_2, \dots)$ is a sequence of unit vectors we let
$\omega_f := \otimes_i \omega_{f_i}$ be the corresponding pure product state of $\F$.
Thus
\[
\omega_f(v_{s_1}\dotsm v_{s_k}v_{t_k}^*\dotsm v_{t_1}^*)
= \ip{v_{s_1}}{f_1}\dotsm\ip{v_{s_k}}{f_k}\ip{f_k}{v_{t_k}}\dotsm\ip{f_1}{v_{t_1}}.
\]

\medskip

\subsection*{A. Periodic pure essential product states}

Suppose $\omega_f$ has finite period $p$; this is equivalent to $p$ being the smallest
positive integer for which the series $\sum(1-\abs{\ip{f_i}{f_{i+p}}})$ converges
\cite[\S4]{laca1}.
The GNS triple for $\omega_f$ is unitarily equivalent to $(\pi_0',\Hh_0',\xi_0')$,
where $\Hh_0'$ is the infinite tensor product $\Ee^{\otimes\infty}$ with canonical
unit vector $\xi_0' := f_1\otimes f_2 \otimes\dots$, and
\begin{multline*}
\pi_0'(v_{s_1}\dots v_{s_k}v_{t_k}^*\dots v_{t_1}^*)(h_1\otimes h_2 \otimes \dots) \\
= \ip{h_1}{v_{t_1}}\dots\ip{h_k}{v_{t_k}}v_{s_1}\otimes\dots\otimes v_{s_k}
\otimes h_{k+1}\otimes h_{k+2} \otimes \dots.
\end{multline*}
(See \cite{von,gui} for the definitions and basic properties of 
infinite tensor products.)
The state $\beta^{*i}\omega_f$ corresponds to the sequence
$(\underbrace{v_1,\dots,v_1}_i,f_1,f_2,\dots)$,
so we can similarly define $(\pi_i',\Hh_i',\xi_i')$.
Replacing $(\pi_i,\Hh_i,\xi_i)$ with $(\pi_i',\Hh_i',\xi_i')$ for $0 \le i \le p-1$
in Theorem~\ref{theorem:gns},
\[
\xi_p' := \underbrace{v_1\otimes\dots\otimes v_1}_p\otimes f_1 \otimes f_2 \otimes\dots
\in\Hh_0'
\]
is a linking vector for $\omega_f$.
For this choice of $\xi_p'$, the vectors $\xi_k$ for $k\ge p+1$
are similarly given by
\[
\xi_k' := \underbrace{v_1\otimes\dots\otimes v_1}_k\otimes f_1 \otimes f_2 \otimes\dots
\in\Hh_k'.
\]
It is routine to check that \eqref{eq:tilderho} yields the same formula
for extensions of $\omega_f$ as that given in \cite[Theorem~3.1]{fowler}.

\medskip

\subsection*{B. Generalized Cuntz states.}
Next we use periodic sequences to construct and classify examples along the lines of 
those from \cite[\S4]{laca1} and \cite[Corollary 5.5]{laca2}.
When $f$ is a constant sequence,
the pure extensions of $\omega_f$ to $\On$ are the Cuntz states \cite{cun},
and lead to shifts which admit a pure normal invariant state \cite{powers}.
We consider here the more general case where $f$ has period $p$, so that it is
determined by  the $p$\ndash tuple $(f_1, \dots, f_p)$.
For such a sequence, $\omega_f$ is periodic in a stronger sense than that 
of Definition~\ref{definition:period}: indeed $\alpha^{*p} \omega_f = \omega_f$. 
 
Although the pure extensions of $\omega_f$ are mutually disjoint,
by Corollary~\ref{cor:classify ergodic}(1)
they induce the same endomorphism of $\bh$ up to conjugacy.
In order to compare the endomorphisms coming from two
different sequences we use Corollary~\ref{cor:classify ergodic}(2). 
The criterion is particularly easy to apply in this strictly periodic
situation because the quasi-orbit of $\omega_f$ is obtained by simply 
taking the cyclic permutations of the $p$\ndash tuple $(f_1, \dots, f_p)$.

\begin{cor}
Suppose $f$ and $g$ are periodic sequences of unit vectors in
$n$\ndash dimensional Hilbert space $\Ee$. 
Let $\alpha_f$  (resp. $\alpha_g$) be the ergodic endomorphism  associated to
a pure extension of $\omega_f $ (resp. $\omega_g$).
Then  $\alpha_f$ is conjugate to $\alpha_g$ if and only if
there is a unitary $W$ on $\Ee$ such that the $p$\ndash tuple of
$1$\ndash dimensional subspaces $([W f_1], \dots, [Wf_p])$
is a cyclic permutation of $([g_1], \dots, [g_p])$.
\end{cor}

\begin{proof}
Corollary~\ref{cor:classify ergodic}(2) reduces the question of conjugacy to
finding a quasi-free automorphism $\gamma_W$ of $\F$ that superimposes the quasi-orbit
of $\omega_f$ to that of $\omega_g$.
By \cite[Corollary 5.3]{laca2} the quasi-orbits of $\omega_f \circ \gamma_W$ and
$\omega_g$ coincide if and only if 
the series $\sum_j (1 - \abs{\ip{W f_j}{g_{j+k}}})$ converges for some $k$.
Since the sequences $f$ and $g$ are periodic, this series converges 
if and only if all its terms vanish; i.e. if and only if 
$[Wf_j] = [g_{j+k}]$ for some fixed $k$ and every $j$.
\end{proof}

\begin{remark}
The orbit of the $p$\ndash tuple $([f_1], \dots, [f_p])$
of one-dimensional subspaces of
$\Ee$ under the joint action of cyclic permutations 
and of the unitary group $\Uu(\Ee)$ (acting diagonally on $p$\ndash tuples)
is thus a {\em complete conjugacy invariant\/} for the class of ergodic endomorphisms
arising from pure essential product states of $\F$
which are strictly periodic under $\alpha^*$.

This invariant also classifies the larger class of ergodic endomorphisms 
associated with pure essential product states
which are {\em eventually\/} strictly periodic, 
in the sense that there exists $p \ge 1$ such that for large enough $k$ one
has $\alpha^{*(k+p) }\omega_f = \alpha^{*k}\omega_f$.
\end{remark}

\medskip

\subsection*{C. Pure extensions of diagonal states.} 

Assume $n$ is finite. 
The diagonal $\Dd$ in $\F$ is the abelian subalgebra generated 
by the projections $v_s v_s^*$, where $s$ is any multi-index.
The spectrum $\hat\Dd$ of $\Dd$ is canonically isomorphic to the 
totally disconnected compact space  $\{1, 2, \dots, n\}^{\NN}$.
A {\em rational point\/} in $\hat\Dd$ is one which corresponds to a
sequence which is eventually periodic, and {\em irrational points\/}
correspond to aperiodic sequences \cite{cun}.

When the sequence $f = (f_i)$ consists of basis vectors
(that is, each $f_i \in \{v_k: 1\le k \le n\}$), 
the state $\omega_f$ of $\F$ is
a {\em diagonal\/} pure state; i.e. it corresponds to a point in $\hat \Dd$.
It was observed by Cuntz that if the sequence is aperiodic then 
the state $\omega_f$ has a unique pure extension. 
Using our Corollary~\ref{cor:gns} we can say what happens at the rational points.

\begin{cor}
Suppose $f$ is a sequence of basis vectors with periodic tail,
so that $\omega_f | \Dd$ is a rational point in the spectrum of $\Dd$. 
Then the pure extensions of $\omega_f$ to $\On$ are mutually disjoint 
and indexed by the circle $\TT$ (via the composition of
$e^{2\pi it} \mapsto e^{2\pi it/p}$ and the gauge action).

Setting $n = 2$ gives uncountably many inequivalent
pure extensions of the trace on the Choi
subalgebra of $\Oo_2$ arising from each rational point in $\hat \Dd$. 
\end{cor}

\begin{proof}
The first assertion follows from Corollary~\ref{cor:gns};
the second one is immediate because if a state of $\Oo_2$
restricts to a diagonal state on $\F$, 
then it extends the trace on the Choi algebra \cite{altw,sta-gi}.
\end{proof}

\end{document}